\documentclass[graybox]{svmult}

\usepackage{mathptmx}       
\usepackage{helvet}         
\usepackage{courier}        
\usepackage{type1cm}        
\usepackage{makeidx}         
\usepackage{graphicx}        
\usepackage{multicol}        
\usepackage[bottom]{footmisc}

\usepackage{amssymb}
\usepackage{amsmath}
\usepackage{amscd}
\usepackage{diagbox}
\usepackage{url}

\def\E{\operatorname{E}}
\def\supp{\operatorname{supp}}
\def\dxb{\underline x_b}
\def\gxb{\overline x_b}
\def\da{\underline a}
\def\ga{\overline a}

\makeindex             

\begin{document}

\title*{Transformation semigroups and their  applications}
\titlerunning{Transformation semigroups and their  applications}
\author{Katarzyna Pich{\'o}r and Ryszard Rudnicki}
\institute{
Katarzyna Pich{\'o}r \at
Institute of Mathematics,
University of Silesia in Katowice, Bankowa 14, 40-007 Kato\-wi\-ce, Poland,
\email{katarzyna.pichor@us.edu.pl}
\and
Ryszard Rudnicki \at Institute of Mathematics,
Polish Academy of Sciences, Bankowa 14, 40-007 Katowice, Poland,
\email{ryszard.rudnicki@us.edu.pl}
}
\maketitle

\abstract{In this chapter we present transformation semigroups and their applications.
We begin with Klein's approach to geometry based on  invariants of transformation groups.
Then we present symmetry groups in chemistry and in classical mechanics.
Next we introduce one-parameter semigroups of transformations and their applications in ergodic theory.
Our main subject are one-parameter semigroups of operators, in particular stochastic semigroups.
We present general results on their existence and long-time behaviour.
We also give examples of one-parameter semigroups related to Markov chains, diffusion and processes with jumps.
We focus on the applications of semigroups of operators in biology.
Among other things, we study models of: DNA evolution; growth of erythrocyte population; gene expression; cell cycle;
the movement of bacteria and insects. We also consider models with stochastic noise and different population models.}


\section{Introduction}
\label{s:intro} 
Transformation semigroups are used today in many scientific fields because they allow 
a unified description of various natural processes and convenient formulation of theorems concerning their properties.     
The intensive study of semigroups began in the early XX century as a generalization of the notion of a group. 
The concept of a group was introduced by \'Evariste Galois\index{Galois} in the 1830s to study polynomial equations.
The abstract definition of a group was formulated by Arthur Cayley \cite{Cayley}\index{Cayley} in 1854 but, as in the case of Galois, it referred only to finite groups.
In 1872 Felix Klein\index{Klein} in his famous the \textit{Erlangen program}\index{Erlangen program} 
presented a novel approach to geometry based on the invariants of transformation groups (Sec. \ref{ss:transformation groups in geometry}).

Transformation groups are used  to study symmetries of some geometrical or chemical objects.
In Sec.~\ref{ss:chemistry} we present a short introduction to molecular and crystal symmetries.
The idea of using transformation groups has been successfully applied to describe Newtonian mechanics~\cite{Arnold}. In Sec.~\ref{ss:mech-Newt} we present 
the group of Galilean transformations  which, 
in combination with Newton's principle of determinacy, allow the derivation of equations of motion in inertial systems; and 
the fundamental laws of classical mechanics.

One-parameter semigroups of transformations of some space $X$ are called dynamical systems.
Different types of dynamical systems are studied  depending on the structure of the space $X$.  
For example if $X$ is a topological space, then we study  continuous dynamical systems on $X$; 
if $X$ is a differentiable manifold --- smooth dynamical systems.
In ergodic theory measure preserving dynamical systems are considered. Such systems consist of measurable transformations on a probability space
$(X,\Sigma,\mu)$ invariant with respect to the measure $\mu$~\cite{CFS, Krengel,LiM,PM}. 
Properties of dynamical systems  can be studied using their representations on the space $L^1(X,\Sigma,\mu)$.
In Sec.~\ref{ss:ergodic} we introduce a semigroup of Frobenius--Perron operators on $L^1$
and we show the relations between ergodic properties of a dynamical system and asymptotic behavior of the 
semigroup of Frobenius--Perron operators.

Among the transformation semigroups, the one-parameter semigroups of operators play an important role, 
as they are used to describe and study processes that change over time. 
They are usually associated with evolution equations of the form $x'(t)=Ax(t)$, where $A$ is some linear (often unbounded) operator
 in a Banach space.  We introduce 
in Sec.~\ref{ss:semigroups-operators}  the notion of  $C_0$-semigroups of operators and we show their relation with evolution equations.

Our main subject are stochastic semigroups.
Stochastic semigroups  are semigroups of positive linear operators defined on the space $L^1$  preserving the set of densities.
These semigroups  are generated by partial differential equations of different types and describe the behaviour of the distributions of Markov processes like diffusion processes, piecewise deterministic Markov processes (PDMPs) and hybrid stochastic processes~\cite{Rudnicki-LN}. They are also used to study ergodic properties of dynamical systems~\cite{LiM}.
In Sec.~\ref{ss:s-s semigroups} we introduce the notions of
stochastic and substochastic semigroups and we state general theorems on their existence. In the following sections, we present basic examples of stochastic semigroups related to: Markov chains (Sec.~\ref{ss:Markov-chains}),
diffusion (Sec.~\ref{ss:d-s}), processes with jumps (Sec.~\ref{ss:proc-jump}),
 processes with random switching (Sec.~\ref{ss:diff-with-switching}),
and to stochastic billiards  (Sec.~\ref{ss:proc-non-random-jumps}).

Long-time behaviour of  semigroups of operators is studied in Part~\ref{s:long-time-beh}. 
We begin with 
a continuous version of the classical Perron--Frobenius theorem (Sec.~\ref{ss:gen-case}).
This theorem concerns semigroups of linear operators on the space $\mathbb R^n$.
The study of such semigroups is a good starting point for introducing asymptotic properties of semigroups in Banach spaces 
and provides simple examples.
One of the tools in studying the behaviour of a 
quasi-compact semigroup of operators when 
$t\to\infty$ is its spectral decomposition (see Sec.~\ref{ss:expon-grow}).
In particular, if the semigroup $\{T(t)\}_{t\ge 0}$ 
on a Banach lattice $E$ is also positive and irreducible, then 
there exist $r\in\mathbb R$,  nonzero $x^*\in E$, and  $\alpha\in E^*$ 
such that 
 \begin{equation}
\label{c-d3-int}
\lim_{t\to\infty}e^{-r\,t}T(t)x=x^*\alpha(x)\quad\textrm{for $x\in E$}.
\end{equation}

Asymptotic properties of stochastic semigroups are studied in the next four sections.
We begin with a theorem concerning  uniformly continuous stochastic semigroups
and with the lower function theorem of Lasota and Yorke~\cite{LY} 
(see Sec.~\ref{ss:stoch-semig-s}).
The main results concern stochastic semigroups  bounded from below by integral operators called partially integral semigroups (see Sec.~\ref{ss:partially-int}).
If such a semigroup $\{P(t)\}_{t\ge 0}$ has a unique 
and positive density $f^*$, then it is asymptotically stable, i.e. 
$P(t)f\to f^*$ for each density $f$.
Strengthening slightly the condition for partial integrability, we obtain a theorem about the decomposition of a semigroup into asymptotically stable parts and a sweeping part (see Sec.~\ref{ss:r-as-decomp}).
This theorem and its conclusions are convenient for the study of many advanced stochastic semigroups.
We show in Sec.~\ref{ss:constrictive} that if a stochastic semigroup has some compact properties then it has also an asymptotic decomposition.

Part~\ref{s:biol-appl} is devoted to applications in biological models.
Our aim is to demonstrate both the wealth of applications of semigroups in biological models as well as the usefulness of the methods presented for their study. 
Among other things, we study models of: DNA evolution; growth of erythrocyte population; population growth with stochastic noise; gene expression;
the movement of bacteria and insects.
We also study structured population models 
and more advanced models based on semi-Markov processes --
kangaroo movement and a cell cycle model~\cite{PR-cell-2022}.

\section{Transformation semigroups everywhere}
\label{s:semigroups-everywhere}
\subsection{Basic definitions and examples}
\label{ss:definition}
A \textit{semigroup}\index{semigroup} $(\mathbf G,\cdot)$ is an arbitrary set $\mathbf G$ together with a binary operation $\cdot$ on $\mathbf G$
with the property of associativity, i.e. $(x\cdot y)\cdot z = x\cdot(y\cdot z)$ for all $x,y,z\in\mathbf G$.
An element $e\in\mathbf G$ is called
\textit{identity}\index{identity} if $e\cdot x=x\cdot e=x$ for all $x\in\mathbf G$. From the property of associativity    
it follows that a semigroup has at most one identity element. 
A \textit{monoid}\index{monoid} is a semigroup with an identity.
An element $x$ of the monoid $(\mathbf G,\cdot,e)$ has two-sided inverse $x^{-1}\in\mathbf G$, shortly called \textit{inverse},\index{inverse element} if
$x\cdot x^{-1}=x^{-1}\cdot x=e$.    
A monoid is called a \textit{group}\index{group} if each element of $\mathbf G$ has an inverse.    
A semigroup $(\mathbf G,\cdot)$ is called \textit{commutative}\index{commutative semigroup} if $x\cdot y=y\cdot x$ for all $x,y\in\mathbf G$.

A semigroup of elements
acting as transformations of some space $X$ is called a \textit{transformation semigroup}\index{transformation semigroup} or a 
\textit{semigroup action on a set}.\index{semigroup!action on a set} 
In this case the binary operation on $\mathbf G$ is the composition $\circ$ of transformations
and the identity element is the identity transformation $I$. 

A \textit{topological semigroup}\index{topological!semigroup} $(\mathbf G,\cdot)$ is
a topological space that is also a semigroup such that the semigroup operation $(x,y)\mapsto x\cdot y$ is continuous.
If $(\mathbf G,\cdot)$ is a group, then  it is called a \textit{topological group}\index{topological!group}
if we additionally assume that the inversion map $x\mapsto x^{-1}$ is continuous.

\textit{One-parameter semigroup} of transformations is a family $\mathbf G=\{g_t\}_{t\in T}$ of transformations of a set $X$
such that $T$ is one of the sets $\mathbb N=\{0,1,\dots\}$, $\mathbb Z$, $\mathbb R$, $\mathbb R_+=[0,\infty)$;
$g_s\circ g_t=g_{s+t}$ for $s,t\in T$ and $g_0=I$. If $X$ is a topological space and the map $\varphi\colon T\times X\to X$ given by
$\varphi(t,x)=g_t(x)$ is continuous, then family  $\mathbf G$ is called a $C_0$-\textit{semigroup}.\index{$C_0$-semigroup}
If $T=\mathbb N$ then there is a function $g$ such that $g_t$ is the $t$-th iterate of $g$
and if $T=\mathbb Z$ then additionally $g_{-t}=g_t^{-1}$ for $\,t>0$. If $T=\mathbb N$ and $g\colon X\to X$ is a continuous function, then
$\mathbf G$ is a $C_0$-semigroup; if $T=\mathbb Z$ and $g$ is a homeomorphism of $X$, then $\mathbf  G$ is a $C_0$-group.
One-parameter semigroup of transformations is also called a \textit{dynamical system}. In the theory of dynamical systems  
a one-parameter semigroup is called: \textit{endomorphism}\index{endomorphism} if $T=\mathbb N$; \textit{automorphism}\index{automorphism} if $T=\mathbb Z$;
\textit{flow}\index{flow} if $T=\mathbb R$; and \textit{semiflow}\index{semiflow} if $T=\mathbb R_+$.   
If $X$ has some additional structure then a dynamical system should be compatible with this structure, e.g. if $X$ is equipped with a topology, then 
$\mathbf G$ should be a $C_0$-semigroup.

Transformation semigroups provide a rich source of many examples and counter-examples of the concepts introduced earlier.
The basic examples are linear maps of $\mathbb R^n$ or $\mathbb C^n$.
Such a map is of the form $x\mapsto Ax$ for some $A\in {\rm M}^n(K)$, where  ${\rm M}^n(K)$ denotes the set of all 
$n\times n$ matrices over the field $K=\mathbb R$ or $K=\mathbb C$, respectively.
Since the composition of linear maps corresponds to the matrix multiplication, 
the set ${\rm M}^n(K)$ with operation of multiplication and with identity matrix is a monoid.
Since ${\rm M}^n(K)$ is isomorphic to $K^{n^2}$ and multiplication is a continuous function,
${\rm M}^n(K)$ is a topological semigroup. 
The subset of ${\rm M}^n(K)$ that consists of all invertible matrices denoted by 
${\rm GL}(n,K)$  is a topological group because the operation $A\mapsto A^{-1}$,  $A\in {\rm GL}(n,K)$, is continuous.  
By ${\rm O}(n)$ we denote the subset of ${\rm GL}(n,\mathbb R)$ that consists of all 
\textit{orthogonal matrices},\index{orthogonal matrix} i.e. $A$ is orthogonal if $A^{\mathrm T}=A^{-1}$, where $A^{\mathrm T}$ is the transpose of $A$.
The set ${\rm O}(n)$ is also a topological group.
Transformations corresponding to orthogonal matrices are linear isometries of $\mathbb R^n$. 

The simplest example of a $C_0$-semigroup is a translation semigroup  on $\mathbb R$ or    $\mathbb R_+$. 
The transformations are given by the formula $g_t(x)=x+t$, where $t\in\mathbb R$ if $X=\mathbb R$
and $t\ge 0$ if $X=\mathbb R_+$. Observe that if $X=\mathbb R$, then the translation semigroup $\{g_t\}$ is a flow and if $X=\mathbb R_+$ it is a semiflow.
Another useful example of a $C_0$-semigroup is the rotation group on a unit circle $S^1$ on a complex space $\mathbb C$.
This group consists of all transformations of the form $g_t(z)=ze^{it}$, $t\in\mathbb R$. 

More general example of a $C_0$-semigroup can be defined by using systems of ordinary differential equations. 
Consider a moving particle in an open set 
 $X\subset \mathbb R^d$.
We assume that if a particle is at point $x$ then its velocity is $b(x)$, where $b\colon X\to \mathbb R^d$ is a differentiable function.
It means that if $x(t)$ is its position at time $t$ then the function $x(t)$ satisfies  the following equation
\begin{equation}
\label{M1}
x' (t) = b(x (t)).
\end{equation}
We also assume that the particle does not leave the set $X$.
For each $\bar x\in X$ we denote by
$\pi_t \bar x$ the solution $x(t)$ of 
(\ref{M1}) with the initial condition $x(0)=\bar x$.
Then $\{\pi_t\}_{t\ge 0}$ is a semiflow on $X$. If additionally   
the solution $x(t)$ exists and $x(t)\in X$
for all $t\in \mathbb R$, then $\{\pi_t\}_{t\in \mathbb R}$ is a flow on $X$.

In Sec.~\ref{ss:ergodic} we introduce measure-preserving dynamical systems
which are one-parameter semigroups but may not be $C_0$-semigroups. Such systems play an essential role in ergodic theory.

An important role in the study of properties of algebraic structures plays the notion of homomorphism.
If $(\mathbf G_1,\cdot)$  and $(\mathbf G_2,\circ)$ are semigroups then a \textit{homomorphism}\index{homomorphism} 
from $\mathbf G_1$ to $\mathbf G_2$ is a map $h\colon \mathbf G_1\to \mathbf G_2$
such that 
\begin{equation}
\label{homomorphism1}
h(g_1\cdot g_2)=h(g_1)\circ h(g_2)\quad \text{for all $\,g_1,g_2\in \mathbf G_1$}.
\end{equation}
If $\mathbf G_1$ and $\mathbf G_2$ are monoids, then we additionally assume that a homomorphism maps the identity element of 
 $\mathbf G_1$ to  the identity element of  $\mathbf G_2$.
If $\mathbf G_1$ and $\mathbf G_2$ are groups and $h$ satisfies \eqref{homomorphism1}, then 
the homomorphism $h$ maps the identity to the identity and $h(g^{-1})= (h(g))^{-1}$ for all $g\in \mathbf G_1$.

We now give an example of a homomorphism from a semigroup without linear structure to a semigroup defined on a linear space. 
We will use the presented construction in Sec.~\ref{ss:ergodic}.    
Let $X$ be any set and let $(\mathbf G,\circ)$ be the semigroup  which consists of all functions from $X$ to $X$  with composition as the binary operation:
$(\varphi,\psi)\mapsto \varphi\circ\psi$.
Let $V=F(X,\mathbb R)$ be the set of all maps from $X$ to $\mathbb R$ and let $\mathbf L=L(V,\mathbb R)$ be the set of all linear operators on the space
$V$. The space $\mathbf L$ is a linear space but it is also a semigroup with the binary operation: $(P,Q)\mapsto Q\circ P$.
Given a function $\varphi\colon X\to X$ we define a linear operator $U_{\varphi}\in \mathbf L$ by the formula
$(U_{\varphi}f)(x)=f(\varphi(x))$ for $f\in V$. Then
\[
U_{\varphi\circ \psi}=U_{\psi}\circ U_{\varphi},
\]
which means that the map $h\colon \mathbf G\to  \mathbf L$ given by $h(\varphi)=U_{\varphi}$ is a homomorphism.
The homomorphism $h$ is a representation of the group $\mathbf G$.
Generally, if $(\mathbf G,\cdot)$ is an arbitrary semigroup,  $V$ is a linear space and  ${\rm End}(V)$
is the set of all linear transformations (endomorphisms) $V\to V$, and $h\colon \mathbf G\to {\rm End}(V)$ is a homomorphism
then $h$ is called  a \textit{representation of the semigroup}\index{representation!of a semigroup} $\mathbf G$ on the vector space $V$.

The operator $U_{\varphi}$ is called the \textit{Koopman operator} or the 
\textit{composition operator} related to $\varphi$.\index{Koopman operator}\index{composition operator}  
Usually the set $X$ has an additional  structure. If for example $(X,\Sigma)$ is a measurable space, then we can consider a subgroup of $\mathbf G$ consisting of all measurable transformations. Then the space $V$ can be the set of all measurable functions from $X$ to $\mathbb R$.

If instead of the above semigroup $\mathbf G$ we consider the group $({\rm S}(X),\circ)$ of all invertible functions from $X$ onto $X$, then the homomorphism $h$
transfers  ${\rm S}(X)$ into the set of invertible linear transformations ${\rm GL}(V)$, $V=F(X,\mathbb R)$, and the homomorphism $h$ is 
a representation of the group ${\rm S}(X)$. In general, if $(\mathbf G,\cdot)$ is an arbitrary group,  $V$ is a linear space,  ${\rm GL}(V)$
is the set of bijective linear transformations from $V$ to $V$, and $h\colon \mathbf G\to {\rm GL}(V)$ is a homomorphism,
then $h$ is called  a \textit{representation of the group}\index{representation!of a group} $\mathbf G$ on the vector space $V$.

\subsection{Transformation groups in geometry}
\label{ss:transformation groups in geometry}
As we have mentioned in Introduction, different geometries can be classified by  transformation groups.
Consider a given group $\mathbf G$ of transformations of the space $\mathbb R^3$. 
Klein proposed that a geometry related to the  
group $\mathbf G$ consists of the objects (properties or figures) 
which, as a result of the action of any transformation $S\in \mathbf G$,  transfer to objects of the same type.
\textit{Euclidean geometry}\index{Euclidean geometry} corresponds to 
the \textit{isometric transformations}\index{isometric transformation} (translations, rotations, reflections and their compositions).
We recall that any isometric transformation is given by the formula
$Sx=Ax+b$, where $b\in\mathbb R^3$ and the matrix $A$  is orthogonal.
The objects of Euclidean geometry are for example: the distance, perpendicularity and parallelity of straight lines, circles, preservation of angles.    

If we replace the group of isometries by the group of \textit{similarities}\index{similarity}, i.e. the transformations 
of the form $Sx=cAx+b$, where $c\ne 0$ is a real constant and $A$ is an orthogonal matrix,
then we obtain another geometry that Euclid actually dealt with.
This geometry has the same aforementioned objects as Euclidean geometry with the exception of distance.

\textit{Affine geometry}\index{affine!geometry} is the geometry of affine transformations, i.e.
the maps of the form $Sx=Ax+b$, where $A$ is any invertible linear transformation. The objects in affine geometry are 
parallelity of straight lines and ratios of lengths. 

More advanced is the \textit{projective geometry}\index{projective!geometry} which is induced by homographies.
A point in the  \textit{projective space}\index{projective!space} $X$ of dimension $n$ is a one-dimensional subspace of the vector space $\mathbb R^{n+1}$.
If $\mathbf x=[x_0,x_1,\dots,x_n]$ is a non-zero element of $\mathbb R^{n+1}$ then 
\[
[x_0:x_1:\dots:x_n]=\{[\theta x_0,\theta x_1,\dots,\theta x_n]\colon \theta \in\mathbb R\}
\]
is a one-dimensional subspace of $\mathbb R^{n+1}$ spanned by the vector $\mathbf x$.
Let 
\[
\mathbb R^{n+1}_0=\mathbb R^{n+1}\setminus\{\mathbf 0\}, \quad \mathbb R^{n+1}_{0,\bullet}=\mathbb R^{n+1}\setminus\{\mathbf x\in\mathbb R^{n+1} \colon x_0=0\}.
\] 
Then we can consider a map $P\colon \mathbb R^{n+1}_0\to X$ given by
\[
P([x_0,x_1,\dots,x_n])=[x_0:x_1:\dots:x_n].
\]
Let $A\colon \mathbb  R^{n+1}\to \mathbb  R^{n+1}$ be an invertible linear map.
A \textit{homography}\index{homography}, called also \textit{projective transformation},\index{projective!transformation} 
induced by  $A$
 is a mapping $f_A$ from $X$ to $X$ given by
\[
f_A([x_0:x_1:\dots:x_n])=[(A\mathbf x)_0:(A\mathbf x)_1:\dots:(A\mathbf x)_n].
\]
If $[a_{ij}]$ is a matrix of the map $A$, then it is called the \textit{matrix of the homography}\index{matrix of the homography} $f_A$.
The following diagram commutes:
\[
\begin{CD}
\mathbb R^{n+1}_0 @> A  >> \mathbb R^{n+1}_0   \\
@VPVV @VVPV\\
X @>>f_A>X,
\end{CD}
\]
i.e. $f_A\circ P=P\circ A$.
If $V$ is a linear subspace of dimension $m+1\le n$ in $\mathbb  R^{n+1}$,
then $P(V\setminus\{\mathbf 0\})$ is called a \textit{flat}\index{flat} of dimension $m$ in $X$.
In particular the lines in $X$ are images under $P$ 
of two-dimensional linear subspaces.
From the formula $f_A\circ P=P\circ A$ it follows that a homography transfers    
a flat on a flat of the same dimension.

\begin{figure}
\begin{center}
\begin{picture}(180,120)(80,-10)
\includegraphics[scale=1,viewport=211 585 393 450]{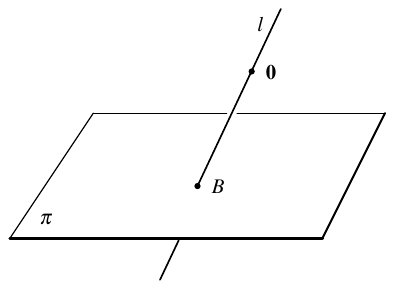}
\end{picture}
\end{center}
\caption{The map $Q$ transfers the line $l$ to the point $B$.}
\label{r:im-1}
\end{figure}

We can also introduce a projective space in a way that is closer to geometrical intuition.
Let $\pi$ be the hyperplane $\{\mathbf x\in \mathbb R^{n+1}\colon x_0=1\}$. Fix a point $\mathbf x\in \mathbb R^{n+1}_{0,\bullet}$. Then 
the vector $\overrightarrow{\,\mathbf x\mathbf 0\,}$ spans the line 
$l=[x_0:x_1:\dots:x_n]$. Denote by $B$  the intersection of $l$ and $\pi$ (see Fig.~\ref{r:im-1}).
In this way we obtain a bijective transformation of the set of all lines which 
are not parallel to the hyperplane $\pi$ onto the space $\mathbb R^n$ (we omit $1$ - the first coordinate of the point $B$).  
The projective space $Y$ is defined by adding the set $Y_{\infty}$ of "points at infinity"  to $\mathbb R^n$, which corresponds to lines  
parallel to the hyperplane $\pi$. Thus $Y=\mathbb R^n\cup Y_{\infty}$.
In one-dimensional case $Y_{\infty}$ is a single point, so we can write $Y=\mathbb R\cup\{\infty\}$.
If $n=2$, then $\mathbb R^2$ and $Y_{\infty}$ can be identified with a complex plane $\mathbb C$ and 
with a semicircle at infinity, respectively, so we can write  
$Y_{\infty}=\{\infty\cdot e^{it}\colon t\in [0,\pi)\}$.
For arbitrary $n$ the set $Y_{\infty}$ is an $(n-1)$-dimensional hemisphere at infinity, 
precisely, the hemisphere is formed by identification of antipodal points on the sphere.    

Let $\mathbf x\in \mathbb R^{n+1}_{0,\bullet}$. Then $B=\big(1,\frac{x_1}{x_0},\dots,\frac{x_n}{x_0}\big)$
and we can define the map $Q\colon \mathbb R^{n+1}_{0,\bullet}\to \mathbb R^n$ by the formula
 $Q(x_0,x_1,\dots,x_n)=\big(\frac{x_1}{x_0},\dots,\frac{x_n}{x_0}\big)$.  
It is not difficult to check that if $V$ is an $(m+1)$-dimensional subspace of $\mathbb R^{n+1}$ 
such that 
$V\setminus \{\mathbf 0\}\subset \mathbb R^{n+1}_{0,\bullet}$, then 
$Q(V\setminus \{\mathbf 0\})$ is an $m$-dimensional flat in $\mathbb R^n$, i.e. $Q(V\setminus \{\mathbf 0\})=W+\mathbf b$,
where $W$ is an $m$-dimensional subspace of $\mathbb R^n$ and $\mathbf b\in\mathbb R^n$.

In this model the homography $f_A$ restricted to the set $\mathbb R^n$ is given by the formula
\[
f_A(x_1,\dots,x_n)=(y_1,\dots,y_n), \quad\text{where $y_k=\frac{a_{k0}+a_{k1}x_1+\dots+a_{kn}x_n}{a_{00}+a_{01}x_1+\dots+a_{0n}x_n}$,\quad $k=1,\dots,n$,}
\]
for all points  $(x_1,\dots,x_n)$ such that the last denominator differs from $0$, otherwise the point is transfered on a point at infinity.
The group of all homographies contains affine maps and central perspectivities.
It is interesting that though a homography is generally a nonlinear map, it transfers  a flat on a flat,
but does not preserve parallelity. 
Since the group of homographies is larger than the group of affine maps 
all theorems of projective geometry remain valid in affine geometry.

It should be noted that $n$-dimensional projective space can be treated as a compact smooth manifold $M$
and construct as the closed $n$-dimensional ball with identified antipodal points on the boundary.
The map $h(\mathbf x)=\frac{\mathbf x}{1+\|\mathbf x\|}$ transfers part $\mathbb R^n$ of $Y$ onto the open ball $B(\mathbf 0,1)$
and the infinite points are transfered on the boundary of $M$. When $n=1$ then the set $Y_{\infty}$ is a one-point set, thus 
the projective space is topologically equivalent to a circle. In $n=2$  the manifold $M$ cannot be embedded in~$\mathbb R^3$.

The idea of using transformation groups to define invariants has found applications in other branches of geometry, for example in
elliptic geometry. In topology that has separated from geometry, we consider objects invariant with respect to a group of homeomorphisms. 
Moreover, homotopy groups and cohomotopy groups play important role in algebraic topology.

\subsection{Symmetry groups in chemistry}
\label{ss:chemistry}
Transformation groups are also applied to study symmetries of some geometrical or chemical objects.
Planar or spatial figures are invariant with respect to some transformations called symmetries.
All symmetries of a given figure $X$ form a subgroup  of ${\rm O}(n)$, $n=2,3$, called a \textit{symmetry group}\index{symmetry group}
and denoted by ${\rm Sym}(X)$. For example, if $X$ is a non-square rectangle then this group consists of four symmetries:
identity $I$, two reflections about the horizontal and vertical middle lines  $S_h$, $S_v$ and the central symmetry $S_O$. 
The map $h\colon {\rm Sym}(X)\to \{0,1,2,3\}$ given by the formulae: $h(I)=0$, $h(S_h)=1$, $h(S_v)=3$, $h(S_O)=2$,
is the homomorphism between the groups ${\rm Sym}(X)$ and $(\mathbb Z_4,\oplus)$, where $\mathbb Z_4=\{0,1,2,3\}$ 
and $a \oplus b= (a+b)\mod 4$.

In chemistry symmetry groups describe molecular and crystal symmetries.
Group theory enables simplification of analysis of chemical and physical 
properties of molecules and crystals.
Molecular symmetries are operations performed on chemical molecules leading to an overlapping arrangement of the molecule's atoms before and after the operation.

\begin{figure}
\begin{center}
\begin{picture}(200,140)(10,20)
\includegraphics[scale=1,viewport=190 567 409 433]{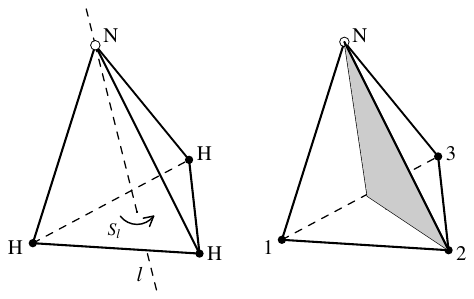}
\end{picture}
\end{center}
\caption{Rotation $S_l$ around the axis $l$ by an angle of $120^{\circ}$ does not change an ammonia molecule. The mirror plane $\sigma_2$ is marked on the right picture.}
\label{r:im-2}
\end{figure}

As an example consider an ammonia molecule ${\rm{NH}}_3$ (see Fig.~\ref{r:im-2}). 
Let $l$ be a straight line passing through the nitrogen atom and the centre of the triangle defined by the hydrogen atoms.  
Then the symmetry group of this molecule consists of six isometries: identity $I$; two rotations $S_1=S_l$, $S_2=S_l^2$ around the axis $l$; and 
three reflections  $S_{\sigma_1}$, $S_{\sigma_2}$, $S_{\sigma_3}$.
Here $S_l$ is the rotation around the axis $l$ by an angle of $120^{\circ}$.
The mirror planes $\sigma_1$, $\sigma_2$, $\sigma_3$ ran through the axis $l$ and one of the hydrogen atoms.  
The symmetry group of ${\rm NH}_3$ is denoted by ${\rm C}_{\rm 3v}$ (see Tab.~\ref{tab-amoniak}). It is obvious that this group is isomorphic to
the symmetry group of the equilateral triangle.

\begin{table}[ht]
\begin{center}
\begin{tabular}{|c|c|c|c|c|c|c|}
\hline
\diagbox[width=3em]{$g$}{$f$}& $I$            &$S_l$              & $S_l^2$          &\,$S_{\sigma_1}$\, &\,$S_{\sigma_2}$\, &\,$S_{\sigma_3}$\,\\
\hline
$I$                     &  $I$                &$S_l$              & $S_l^2$          &$S_{\sigma_1}$ &$S_{\sigma_2}$ &$S_{\sigma_3}$\\[0.5ex]
$S_l$                   &  $S_l$              &$S_l^2$            & $I$              &$S_{\sigma_3}$ &$S_{\sigma_1}$ &$S_{\sigma_2}$\\[0.5ex]
$S_l^2$                 &  $S_l^2$  	      &$I$                & $S_l$            &$S_{\sigma_2}$ &$S_{\sigma_3}$ &$S_{\sigma_1}$\\[0.5ex]
$S_{\sigma_1}$          &  $\,S_{\sigma_1}$\, &\,$S_{\sigma_2}$\, &$\,S_{\sigma_3}$\,& $I$           &$S_l$          &$S_l^2$       \\[0.5ex]
$S_{\sigma_2}$          &  $S_{\sigma_2}$     &$S_{\sigma_3}$     &$S_{\sigma_1}$&$S_l^2$        & $I$           &$S_l$         \\[0.5ex]
$S_{\sigma_3}$          &  $S_{\sigma_3}$     &$S_{\sigma_1}$     &$S_{\sigma_2}$&  $S_l$        &$S_l^2$        & $I$          \\[0.5ex]
\hline
\end{tabular}
\end{center}
\caption{The results of the operation  $g\circ f$ in the group  ${\rm C}_{\rm 3v}$.}
\label{tab-amoniak}
\end{table} 

A group of transformations of the spaces $\mathbb R^n$ is called a \textit{point group}\index{point group} 
if all transformations have a common fixed point.
All isometries of a bounded object are point groups. Thus, a molecular symmetry group is also a point group.
A \textit{crystallographic group}\index{crystallographic!group} 
is the group of all isometries that leave the crystal structure  unchanged,
i.e. the same kinds of atoms would be placed in similar positions as before the transformation.
A crystallographic group  includes some translations of the space.
A \textit{crystallographic point group}\index{crystallographic!point group} is a subgroup  of the crystallographic group consisting of all transformations
having a common fixed point.
Thought there are infinitely many three-dimensional point groups,  there are only 32 crystallographic point groups.

The crystallographic point group determines,  the directional variation of physical properties that arise from its structure, including optical
and electro-optical features. The invariance of a crystal with respect to a translational defines crystallinity.
Other applications of symmetry groups in quantum mechanical and in the spectroscopy  can be found in \cite{Blinder,Carter,Cotton}.

\subsection{Applications to Newtonian mechanics}
\label{ss:mech-Newt}
We will study the motion of material points
in a three-dimensional Euclidean space.  
The position of points changes with time, so it is natural to assume
that time is an element of space. Thus, physical space, called
\textit{world}\index{world}, is a four-dimensional space. Since in this
space we do not distinguish the origin of the coordinate system,
the world is considered to be a four-dimensional \textit{affine space}\index{affine!space} 
 $A^4=A^1 \times A^3$,
in which $A^1$ is identified with the time axis, and at each fixed
 time $t$ a set of points (events) of the form
 $\{ t\} \times A^3$  is called \textit{subspace of 
 simultaneous events}\index{subspace of simultaneous events}.
 The difference of points from the space $A^3$
 is a vector from $ \mathbb{R}^3$, so in the  
subspace of simultaneous events it is possible to introduce a distance
\[
\varrho( ( t, \mathbf a ),  (t, \mathbf b )) = \| \mathbf a -\mathbf b \| =
\sqrt {\smash[b]{(  \mathbf a - \mathbf b, \mathbf a - \mathbf b)}},
\]
where $ \| \cdot \| $ ,  $( \cdot , \cdot )$ 
is respectively the norm and the scalar product in
$\mathbb{R}^3$.

Many coordinate systems can be introduced in the space $A^4$,
among which we distinguish \textit{inertial systems}\index{inertial system}.
These systems are characterised by 
\textit{Galileo's principle of relativity}:\index{Galileo's principle of relativity}
\begin{enumerate}
\item[(1)] all the laws of nature are the same at all moments of times in all
inertial coordinate systems,
\item[(2)] each coordinate system moving in uniform rectilinear motion with respect to an inertial system is also an
inertial system.
\end{enumerate}
We also assume 
\textit{Newton's principle of determinacy}\index{Newton's!principle of determinacy}, which states that the positions and velocities
of the points of the system at the initial moment uniquely determine the entire motion of the system.

 We introduce the coordinate system in the space 
 $\mathbb{R} \times \mathbb{R}^3$.
Then the motion of a single point can be described by the function
$\mathbf{x}\colon I \to \mathbb{R}^3$, 
where $I$ is the interval
on the time axis. 
The set  $\mathbf{x} (I)$
is called the \textit{trajectory}\index{trajectory} of a point,
and the graph of the mapping, i.e. the set
$\{ ( t, \mathbf{x}(t) )\colon t \in I \}$,
is called the \textit{world line}.\index{world!line}  

The following transformations of the
space $\mathbb{R}\times \mathbb{R}^3$ have the property that the world lines of any mechanical system subjected to the same
transformation pass to other world lines of the same system:
\begin{enumerate}
\item [(a)] uniform motion with velocity $\mathbf v \in \mathbb{R}^3$
\[
g_1 (t, \mathbf x ) = ( t,\mathbf x + \mathbf v t ),
\]
\item [(b)] translation of the origin 
$(s,\mathbf a) \in \mathbb{R} \times \mathbb{R}^3$ 
\[
g_2 (t, \mathbf x ) = (t+s, \mathbf x +\mathbf a),
\]
\item[(c)] orthogonal transformation (rotation or plane symmetry)
\[
g_3 (t, \mathbf x) = ( t, G\mathbf x),
\]
\end{enumerate}
where $ G:  \mathbb{R}^3 \to \mathbb{R}^3$
is an orthogonal transformation. 
The set of all possible compositions of transformations of types $g_1, g_2, g_3$ forms the 
\textit{group of Galilean transformations}\index{group!of Galilean transformations}.
It can be easily checked that any Galilean transformation is 
of the form $g_1\circ g_2\circ g_3$. Preserving the world line
by the $g_1$ transformation follows from Galileo's principle of relativity,
while by the transformations $g_2$ and $g_3$ from the geometrical structure of
space: invariance to displacements (homogeneity) 
and rotations (isotropy).

It is possible to derive the general form of the
equations of motion of material points from Newton's principle of determinacy. Since the motion of points is
uniquely determined by their initial positions
and velocities, so the acceleration at time $t$ is a function of the position and
velocity also at time $t$. Consider $n$ material points.
Let $\mathbf x_i (t)$ describes the position of the $i$-th point at time $t$.
Then the function
\[
\mathbf x(t) = ( \mathbf x_1 (t), \ldots, \mathbf x_n (t) )
\]
describes the motion of a system of $n$ points. Let $N=3n$.
Since acceleration is a function of position and velocity,
so there exists a function
$\mathbf F:  \mathbb{R}^N  \times \mathbb{R}^N \times \mathbb{R} \to \mathbb{R}^N$ such that
\begin{equation}
\label{e:MN1}
\mathbf x''(t) = \mathbf F (\mathbf x (t), \mathbf x^{\prime}  (t), t ),
\end{equation}
where $\mathbf x^{\prime} (t)= (\mathbf x_1^{\prime} (t),\ldots, \mathbf x_n^{\prime}(t) )$.

We will now show how the invariance of the world line with respect to the
Galilean group imposes certain conditions on the function $\mathbf F$. One
of the Galilean transformations is a time shift. This means that
if $\mathbf x (t)$ is a solution of Eq.~\eqref{e:MN1}, then also
for any $s$ the function $\mathbf y (t) =\mathbf x (t+s)$ is a
solution of Eq.~\eqref{e:MN1}.
From Eq.~\eqref{e:MN1} we obtain
\begin{equation}
\label{e:MN2}
\mathbf x''(t+s) = \mathbf F (\mathbf x (t+s), \mathbf x^{\prime} (t+s), t+s).
\end{equation}
The function $\mathbf y$ also satisfies Eq.~\eqref{e:MN1},
thus
\begin{equation}
\label{e:MN3}
\mathbf y'' (t) = \mathbf F (\mathbf y (t), \mathbf y' (t), t )=
\mathbf F ( \mathbf x (t+s), \mathbf x^{\prime} (t+s), t).
\end{equation}
Since $\mathbf x'' (t+s) = \mathbf y'' (t)$, from Eqs.~\eqref{e:MN2} and~\eqref{e:MN3} it follows that
\[
 \mathbf F (\mathbf x (t+s), \mathbf x^{\prime} (t+s), t+s) = 
 \mathbf F (\mathbf x (t+s), \mathbf x^{\prime} (t+s), t).
\]
The initial positions and initial velocities of 
points can be chosen in any way and at any
moment, so
\[
\mathbf F ( \mathbf x, \mathbf y, t+s) = F(\mathbf x, \mathbf y, t)
\]
for all $\mathbf x \in \mathbb{R}^N$, $\mathbf y \in \mathbb{R}^N$,
$t\in  \mathbb{R}$ and $s  \in \mathbb{R}$.
Thus the function $\mathbf F$ does not depend on time
and the equation of motion is of the form
\begin{equation}
\label{e:MN4}
\mathbf x''(t)= F (\mathbf x(t),\mathbf x'(t) ).
\end{equation}
Fix now a vector $\mathbf a \in \mathbb{R}^3$. If
$\mathbf x(t)$ satisfies Eq.~\eqref{e:MN1}, then also 
the function $\mathbf y(t) = (\mathbf x_1 (t) + \mathbf a, \ldots,
\mathbf x_n (t) + \mathbf a)$ satisfies Eq.~\eqref{e:MN1}.
Thus, as before, we conclude that
\[
\mathbf F (\mathbf x_1 + \mathbf a,\ldots, \mathbf x_n +\mathbf a, \mathbf y)=
\mathbf F ( \mathbf x_1,\ldots, \mathbf x_n, \mathbf y),
\]
where $\mathbf x_1 \in \mathbb{R}^3,\ldots, \mathbf x_n \in
\mathbb{R}^3$ and $\mathbf y \in \mathbb{R}^N$.
From the invariance of the transformation
$g_1 ( t,\mathbf x) = ( t,\mathbf x +\mathbf v t)$ 
we conclude that
\[
\mathbf F ( \mathbf x, \mathbf y_1+\mathbf v,\ldots, \mathbf y_n + \mathbf v)=
\mathbf F (\mathbf x, \mathbf y_1,\ldots, \mathbf y_n)
\]
for all $\mathbf x \in \mathbb{R}^N$ and
$\mathbf y_1  \in \mathbb{R}^3,\ldots,
\mathbf y_n \in \mathbb{R}^3$. 
If $\mathbf x(t)$ is a
solution of Eq.~\eqref{e:MN1}, then 
from invariance with respect to 
$g_3 (t,\mathbf x) = (t, G \mathbf x)$, 
the function 
$\mathbf y (t) =(G \mathbf x_1 (t),\ldots, G \mathbf x_n (t) )$
is also a solution of  Eq.~\eqref{e:MN1}. Thus
\begin{equation}
\label{e:MN5}
( G \mathbf x_1^{\prime \prime},\ldots, G \mathbf x_n^{\prime\prime})=
\mathbf F  (  G \mathbf x_1,\ldots, G \mathbf x_n,
 G \mathbf x_1^{\prime},\ldots, G \mathbf x_n^{\prime} ).
\end{equation}
Since
 \[ (G \mathbf x_1^{\prime \prime},\ldots, G \mathbf x_n^{\prime\prime} )=
( G \mathbf F_1 ( \mathbf x, \mathbf x'),\ldots,
G \mathbf F_n (\mathbf x, \mathbf x' ) ),
\]
we obtain from
Eq.~\eqref{e:MN5} that
\[
\mathbf F  (  G \mathbf x_1,\ldots, G \mathbf x_n,
 G \mathbf y_1,\ldots, G \mathbf y_n ) =
(G \mathbf F_1 (\mathbf x, \mathbf y),\ldots,G\mathbf F_n (\mathbf x, \mathbf y) )
\]
for all $\mathbf x_1 \in \mathbb{R}^3,\ldots, \mathbf x_n \in   \mathbb{R}^3$, 
$\mathbf y_1\in \mathbb{R}^3,\ldots,\mathbf  y_n \in \mathbb{R}^3$.
Let us introduce the following notations. If
$\mathbf a \in \mathbb{R}^3$ and $\mathbf x =
(\mathbf x_1,\ldots, \mathbf x_n) \in \mathbb{R}^N$, then
$\mathbf x + \mathbf a: = (\mathbf x_1 + \mathbf a,\ldots,\mathbf x_n +\mathbf a)$ and
if $G: \mathbb{R}^3 \to \mathbb{R}^3$, then
$G\mathbf x:= (G \mathbf x_1,\ldots, G \mathbf x_n)$.
Then the function  $\mathbf F$ satisfies the following conditions
\begin{enumerate}
\item[(i)]  
 $\mathbf F (\mathbf x, \mathbf y) = \mathbf F(\mathbf x + \mathbf a,
\mathbf y + \mathbf b)$,
\item[(ii)]
$G \mathbf F( \mathbf x, \mathbf y) = \mathbf F (G\mathbf x, G \mathbf y)$
\end{enumerate}
for all vectors $\mathbf a \in \mathbb{R}^3, \mathbf b \in \mathbb{R}^3$, 
and for all  orthogonal transformations
$G: \mathbb{R}^3 \to \mathbb{R}^3$.

We now give an example of the application of the obtained result.
Let us consider a system consisting of only one point.
Then it follows from equation (i) that there exists a vector 
$\mathbf c \in \mathbb{R}^3$ such that
  $\mathbf F (\mathbf x, \mathbf y) = \mathbf c$ for each
$\mathbf x \in \mathbb{R}^3$ and $\mathbf y \in \mathbb{R}^3$.
From (ii) it follows that 
we have $G \mathbf c = \mathbf c$
 for any rotation  $G$.
Hence  $\mathbf c = \mathbf 0$. 
Thus, a point moves according to the equation
\[
\mathbf x^{\prime\prime} (t)= \mathbf 0.
\]
It implies that $\mathbf x (t) =\mathbf  a t + \mathbf b$ for some vectors  $\mathbf a \in \mathbb{R}^3$
and $\mathbf b \in \mathbb{R}^3$. 
Thus, a point moves in a straight line at constant velocity. We have proved
\textit{Newton's first law of motion}:\index{Newton's!first law of motion}
if no force acts on a material point, then the point remains at rest or keeps moving in a straight line at constant speed.

\subsection{Measure-preserving dynamical systems}
\label{ss:ergodic}
Let $(X,\Sigma,\mu)$ be a probability space, i.e. $\mu$ is a probability measure defined on a $\sigma$-algebra $\Sigma$ of subsets of $X$.
Consider a one-parameter semigroup
$\{\pi_t\}_{t\in T}$.
We assume that the function $(t,x)\mapsto \pi_t x$ is measurable as a function from  the Cartesian product $T\times X$ to $X$
and that the measure $\mu$ is
\textit{invariant}\index{invariant!measure} with respect  to each transformation $\pi_t$, $t\in T$, i.e.
 $\mu(\pi_t^{-1}(A))=\mu(A)$ for any set $A\in\Sigma$.
Then the quadruple $(X,\Sigma,\mu,\pi_t)$ is called
 a \textit{measure-preserving dynamical system}.
\index{measure-preserving dynamical system}
The triple $(X,\Sigma,\mu)$ is called a \textit{phase space}\index{phase space}
and the set
 $\{\pi_tx\colon t\in T\}$ 
is called 
\textit{trajectory of a point}\index{trajectory!of a point}
$x\in X$.
As in Sec.~\ref{ss:definition}, we divide 
measure-preserving dynamical systems into  
endomorphisms, automorphisms, flows, and semiflows.

Our goal is to study properties of measure-preserving dynamical systems:
ergodicity, mixing and exactness. 
A set $A$ is called  \textit{invariant}\index{invariant!set} with respect to a 
semigroup $\{\pi_t\}_{t\in T}$,
if $\pi_t^{-1}(A)=A$ for any $A\in \Sigma$ and $t\in T$. 
The family of invariant sets forms a $\sigma$-algebra  $\Sigma_{\mathrm{inv}}$.\index{Sigmaa@$\Sigma_{\mathrm{inv}}$}
If $\sigma$-algebra $\Sigma_{\mathrm{inv}}$ is 
trivial, i.e. it consists only of sets of measure $\mu$ zero or one, then   
the  measure-preserving dynamical system
 is said to be \textit{ergodic}\index{ergodic dynamical system}.   

A stronger property than ergodicity is mixing. 
A measure-preserving dynamical system  $(X,\Sigma,\mu,\pi_t)$  is called \textit{mixing}\index{mixing} if 
\begin{equation}
\label{d:mixing}
\lim_{t\to\infty} \mu(A\cap \pi_t^{-1}(B))=\mu(A)\mu(B)\quad 
\text{for  all  $A, B \in \Sigma$}.
\end{equation}
Identifying the measure $\mu$ with the probability P,  one can formulate the condition
\eqref{d:mixing} as follows:
\[   
\lim_{t\to\infty} {\rm P}(\pi_t(x)\in A|x\in B)={\rm P}(A)\quad
\text{for  all $A,B \in \Sigma$ and ${\rm P}(B)>0$},
\]
which means that 
the trajectories of almost all points enter 
a set $A$ with  asymptotic probability $\mu(A)$.

A stronger property than mixing is exactness.
A system $(X,\Sigma,\mu,\pi_t)$  with  double measurable transformations $\pi_t$, i.e. $\pi_t(A)\in \Sigma$ and 
$\pi_t^{-1}(A)\in \Sigma$ for all $A\in \Sigma$ and $t\in T$,  and with  an invariant probability measure 
$\mu$ is called 
\textit{exact}\index{exact}
if for every set 
$A \in \Sigma$ with $\mu(A)>0$ 
we have $\lim_{t\to\infty} \mu(\pi_t(A))=1$.
Exactness cannot hold for flows and automorphisms because in these cases
\[
\mu(\pi_t(A))=\mu(\pi_t^{-1}(\pi_t(A)))=\mu(A).
\]
Exactness is equivalent to the following condition:
 $\sigma$-algebra 
$\bigcap_{t \ge 0}  \pi_t^{-1}(\Sigma)$ 
contains only sets of measure $\mu$ zero or one.
Here $\pi_t^{-1}(\Sigma)=\{\pi_t^{-1}(A)\colon A\in \Sigma\}$.
\vskip1mm

Directly checking whether a dynamical system is ergodic, mixing or exact is not an easy issue. 
This problem can be reduced to the study of linear operators related to the representation of semigroups. 
In order to do it we introduce the notion of Frobenius--Perron operators.

Let $(X,\Sigma,m)$ be a $\sigma$-finite measure space.
A measurable map $\varphi\colon X\to X$  
is called \textit{non-singular}\index{non-singular map}
if it satisfies the following condition 
\begin{equation}
\label{trana2}
m(A)=0 \Longrightarrow  m(\varphi^{-1}(A))=0\textrm{ for $A\in \Sigma$.}
\end{equation}
Let $L^1=L^1(X,\Sigma,m)$ and let $\varphi$ be a measurable nonsingular transformation of $X$.
An operator  $P_{\varphi}\colon  L^1\to L^1$
which satisfies the following condition
\begin{equation}
\label{def-FP}
\int_A P_{\varphi} f(x)\,m(dx)=
\int_{\varphi^{-1}(A)} f(x)\,m(dx)
\textrm{\,  for $A\in \Sigma$ and $f\in L^1$}
\end{equation}
is called the
{\it Frobenius--Perron operator} for the transformation~$\varphi$. 
The operator $P_{\varphi}$ is linear, \textit{positive}\index{positive!operator} 
(if $f\ge 0$ then $P_{\varphi}f \ge 0$) and preserves the integral    
($\int_XP_{\varphi}f\,dm=\int_Xf\,dm $).
Any operator having the above properties is called 
\textit{stochastic ({\rm or} Markov) operator.}\index{stochastic!operator}\index{Markov!operator}
We can equivalently define a stochastic operator as a linear operator
on $L^1$ carrying out densities on densities.
We recall that $f\in L^1$ is a \textit{density}\index{density}
if $f\ge 0$ and $\|f\|=1$. By $D$ we denote the set of all densities. 

The adjoint of the Frobenius--Perron operator
$P^*\colon L^{\infty}\to L^{\infty}$ is 
given by $P_{\varphi}^*g(x)= g(\varphi(x))$, thus $P^*_{\varphi}$
is the Koopman operator
restricted to the space $L^{\infty}=L^{\infty}(X,\Sigma,m)$.

\begin{example}
\label{e:F-P-Rd}
Let  $X$ be a subset  of  $\mathbb R^d$
with nonempty interior and with the boundary of zero Lebesgue measure. 
Let  $\varphi\colon X\to X$ be a measurable transformation.
We assume that there exists 
at most countable family of pairwise disjoint open sets
$U_i$, $i\in I$, contained in $X$ having the following properties:
\begin{enumerate}
\item[a)] the sets  $X_0=X\setminus \bigcup_{i\in I} U_i$ and
 $\varphi(X_0)$ have zero Lebesgue measure,      
 \item[b)] maps $\varphi_i =\varphi\Big |_{U_i}$ are diffeomorphisms
from $U_i$ onto $\varphi(U_i)$, i.e., 
$\varphi_i$ are $C^1$ and invertible transformations  and 
$\det \varphi'_i(x)\ne 0$
at each point  $x\in U_i$.      
\end{enumerate}

\noindent Then transformations  $\psi_i= \varphi_i^{-1}$ are also 
diffeomorphisms  from $\varphi(U_i)$ onto  $U_i$ and 
the Frobenius--Perron operator $P_{\varphi}$ exists and is given by the formula
\begin{equation}
\label{F-P-operator-gladki}
P_{\varphi}f(x)=\sum_{i\in I_x} f(\psi_i(x))|\det \psi_i'(x)|,
\end{equation} 
where $I_x=\{i\colon x\in \varphi(U_i)\}$ (see \cite{Rudnicki-LN,RT-K-k}). 
\end{example}

\begin{example}
\label{e:F-P-cw}
If we consider a flow $\{\pi_t\}_{t\in \mathbb R}$ on $X$ related to
Eq.~(\ref{M1}),
then according to formula (\ref{F-P-operator-gladki})  the group $\{P^t\}_{t\in \mathbb R}$ of Frobenius--Perron operators is given by the formula
\[
P^tf(x)=f(\pi_{-t}x) \det\Big[ \dfrac{d}{d x} \pi_{-t}x\Big].
\]
\end{example}

The set of all nonsingular transformations with composition as the binary
operation forms a monoid $(\mathbf G,\circ)$. 
The set of all stochastic operators ${\rm StO}(X)$ is also a monoid
with composition as  the binary operation and ${\rm StO}(X)$ is a linear space.
Since
\[
P_{\varphi\circ \psi}=P_{\varphi}\circ P_{\psi},
\]
the map $h\colon \mathbf G\to  {\rm StO}(X)$ given by $h(\varphi)=P_{\varphi}$
 is a representation of the semigroup $\mathbf G$ on the space $L^1(X,\Sigma,m)$.

Frobenius--Perron operators can be successfully used to study
ergodic properties of  dynamical systems \cite{LiM}. 
Let    $\varphi\colon X\to X$ be a nonsingular transformation of some 
$\sigma$-finite measure space $(X,\Sigma,m)$.
First, we check when a probability measure $\mu$ 
 absolutely continuous with respect to $m$ is invariant with respect to $\varphi$. 
Let $f^*$ be the density of $\mu$ with respect to $m$, i.e. 
$\mu(A)=\int_Af^*\,dm$ for $A\in \Sigma$.  
The measure $\mu$ is invariant with respect to $\varphi$  
if and only if 
\[
\int_{\varphi^{-1}(A)} f^*(x)\,m(dx)=\int_A  f^*(x)\,m(dx)
\textrm{\,  for $A\in \Sigma$.}
\] 
From 
\eqref{def-FP}
it follows immediately that the measure $\mu$ is invariant with
respect to $\varphi$  if and only if $P_{\varphi}f^*=f^*$.
It means that if the map
$\pi\colon T\times X\to X$ is measurable and
 $\{\pi_t\}_{t\in T}$ is
a one-parameter semigroup of nonsingular transformations of $(X,\Sigma,m)$,
$P^t$ denotes the Frobenius--Perron operator corresponding to $\pi_t$,
 then the quadruple $(X,\Sigma,\mu,\pi_t)$ 
is a measure-preserving dynamical system
if and only if $P^tf^*=f^*$ for all $t\in T$.

It is also not difficult to check that the measure
$\mu$ is ergodic if and only if  $f^*$ is a unique fixed point of all operators $P^t$
 in the set of densities; and
the dynamical system is mixing if and only if for every $f\in D$ the density  $f^*$ is the weak limit of $P^tf$ as $t\to\infty$.
By the \textit{weak limit}\index{weak limit} $\lim_{t\to\infty} P^tf$ we understand a function $h\in L^1$  such that for 
every $g\in L^{\infty}$ we have 
\[
\lim_{t\to\infty}\int_X  P^tf(x)g(x)\,m(dx)=\int_X  h(x)g(x)\,m(dx)
\]
and we denote it by w-$\lim_{t\to\infty} P^tf$.
More difficult is to prove that
the dynamical system 
is exact if and only if  $\lim_{t\to\infty} P^tf=f^*$ for every $f\in D$. 
 
We collect the relations between ergodic properties of a dynamical system  and the behavior of the 
Frobenius--Perron operator in Table~\ref{tab1}.
\begin{table}[ht]
\begin{center}
\begin{tabular}{|c|c|}
\hline
$\mu$ & $f^*$
\\
\hline
invariant&$P^tf^*=f^*$ for all $t\in T$\\
ergodic&$f^*$ is a unique fixed point  in $D$ of all $P^t$ \\
mixing& w-$\lim_{t\to\infty} P^tf=f^*$ for every $f\in D$\\
exact & $\lim_{t\to\infty} P^tf=f^*$ for every $f\in D$\\
\hline
\end{tabular}
\end{center}
\caption{The relations between ergodic properties of the dynamical system $(X,\Sigma,\mu,\pi_t)$ and  the  Frobenius--Perron operators $P^t$.}
\label{tab1}
\end{table} 

\begin{example}
\label{Tent-map}
Let $\varphi\colon [0,1]\to [0,1]$ be the transformation given by 
\begin{equation}
\label{tent}
\varphi(x)=
\begin{cases}
2x &\textrm{for $x\in [0,1/2]$},\\
2-2x&\textrm{for $x\in (1/2,1]$.}
\end{cases}
\end{equation} 
 The transformation  $\varphi$ is called 
the \textit{tent map}\index{tent map} because its graph looks like a tent.
The Frobenius--Perron operator $P_{\varphi}$ is of the form
\begin{equation}
\label{FP-diad}
P_\varphi f(x)= \tfrac 12 f(\tfrac 12)+\tfrac 12f(1-\tfrac 12x).
\end{equation} 
It is easy to see that the density $f^*=1_{[0,1]}$ satisfies $P_\varphi f^*=f^*$,
and consequently  $\varphi$ is an endomorphism on the phase space $([0,1],\mathcal B([0,1]),\ell)$,
where $\mathcal B([0,1])$ is the $\sigma$-algebra of Borel subsets of $[0,1]$ and $\ell$ is the Lebesgue measure.
We check that $\lim_{t\to\infty}P_{\varphi}^tf=f^*$ for any density $f$.
It is sufficient to check  this condition 
for densities which are Lipschitz continuous. 
Let $L$ be  the Lipschitz constant for $f$. Then
\[
|P_{\varphi} f(x)-P_{\varphi}f(y)|\le \tfrac 12
 |f(\tfrac{x}2) -f(\tfrac{y}2)|   +\tfrac 12  | f(1-\tfrac 12x)+ f(1-\tfrac 12y)|\le  \tfrac L2 |x-y|.   
\]
Thus   $L/2$ is the  Lipschitz constant for  $P_{\varphi}f$ and 
 by induction we conclude that $L/2^t$ is  the  Lipschitz constant for  
$P_{\varphi}^tf$.
Hence, the sequence $(P_{\varphi}^tf)$ converges uniformly to a constant function.
Since $P_{\varphi}^tf$  are densities, $(P_{\varphi}^tf)$ converges to $f^*$ uniformly, which implies 
the convergence in $L^1$.
The condition $\lim_{t\to\infty}P_{\varphi}^tf=f^*$
implies the exactness of the endomorphism $\varphi$.
\end{example}
\begin{example}
\label{log-map}
Now we check that the \textit{logistic map}\index{logistic map}  $\psi(x)=4x(1-x)$ on $[0,1]$ is exact. 
Let $\varphi$ be the tent map and $\alpha(x)=\tfrac12-\tfrac 12\cos(\pi x)$. Then $\psi\circ \alpha=\alpha\circ \varphi$. 
This implies that 
$P_{\psi}P_{\alpha}=P_{\alpha}P_{\varphi}$. Hence
$P_\psi=P_{\alpha}P_{\varphi}P_{\alpha^{-1}}$
and, by induction, $P^t_{\psi}=P_{\alpha}P^t_{\varphi}P_{\alpha^{-1}}$. 
Let $f\in D$. Then     $P_{\alpha^{-1}}f\in D$ and from the previous example
$\lim_{t\to\infty} P^t_{\varphi}P_{\alpha^{-1}}f=1_{[0,1]}$.
Therefore
\[
\lim_{t\to\infty} P^t_{\psi}f=P_{\alpha}1_{[0,1]}=g^*, \quad \text{where $g^*(x)=\Big[\pi\sqrt{x(1-x)}\,\Big]^{-1}$}. 
\]
Thus the map $\psi$
is an exact endomorphism on the space $([0,1],\mathcal B([0,1]),\mu)$, where $d\mu=g^*(x)\,dx$.
\end{example}

More advanced examples of applications of Frobenius--Perron operators to study ergodic properties of dynamical systems 
can be found in~\cite{BG,LiM,LY}. 

\section{Stochastic semigroups}
\label{s:stochastic-semigroups}
\subsection{Semigroups of operators}
\label{ss:semigroups-operators}
Stochastic semigroups belong to a wider family of semigroups called 
$C_0$-semigroups of operators, so  
we first need to introduce the notion of a $C_0$-semigroup 
and to show its relation with evolution equations.

Let $(E,\|\cdot\|)$ be a Banach space and let    $\{T(t)\}_{t\ge0}$ be a family of linear bounded 
operators on $E$. 
The family $\{T(t)\}_{t\ge0}$ is called 
a \textit{semigroup of bounded operators}\index{semigroup!of  bounded operators}
if it satisfies the following conditions:
\begin{enumerate}
\item[(a)]  $T(0)=I$,  i.e., $T(0)x =x$,
\item[(b)] $T(t+s)=T(t) T(s)\quad \textrm{for}\quad
s,\,t\ge0$.
\end{enumerate}
If additionally
\begin{enumerate}
\item[(c)] for each $x\in E$ the function
$t\mapsto T(t)x$ is continuous,
\end{enumerate}
then $\{T(t)\}_{t\ge0}$ is called
a $C_0$-\textit{semigroup}
or \textit{strongly continuous semigroup}\index{$C_0$-semigroup}\index{strongly continuous semigroup}.
Let $\{T(t)\}_{t\ge0}$ be a $C_0$-semigroup 
and let
 $\mathfrak D(A)$ be the set of such  $x\in E$, that there exists the limit
\begin{equation}
\label{generator}
Ax=\lim_{t\to 0^+}\frac{T(t)x-x}{t}.
\end{equation}
Then the set $\mathfrak D(A)$ is a linear subspace dense in  $E$, 
and $A$ is a linear operator from   
$\mathfrak D(A)$ to $E$. The operator $A$ is called 
the \textit{infinitesimal generator}\index{infinitesimal generator}
(briefly the \textit{generator})\index{generator}
of the semigroup $\{T(t)\}_{t\ge 0}$.
We also say that the operator  $A\colon \mathfrak D(A)\to E$ generates 
the semigroup  $\{T(t)\}_{t\ge 0}$.

The notion of a $C_0$-semigroup and its generator is strictly connected 
with differential equations in  the Banach space $E$.
Let $\{T(t)\}_{t\ge 0}$ be a $C_0$-semigroup on $E$ and let 
$A:\mathfrak D(A)\to E$ be its generator. Then for every  
$x_0\in \mathfrak D(A)$  the function $x\colon [0,\infty)\to E$ defined by the formula $x(t)=T(t)x_0$ is differentiable (in the sense of  Fr\'echet) for $t\ge 0$ 
and satisfies the  equation      
\begin{equation}
\label{r-ew}
x'(t)=Ax(t) \quad \textrm{with initial condition $x(0)=x_0$}.
\end{equation}
An equation of the form  (\ref{r-ew}) is called 
an \textit{evolution equation}.\index{evolution equation}
We also say that equation (\ref{r-ew})  \textit {generates semigroup} $\{T(t)\}_{t\ge 0}$.

In the theory of semigroups of operators, the key issue is when a linear operator is a generator of a $C_0$-semigroup.
If $A$ is a bounded linear operator on $E$, then $A$ generates
a $C_0$-semigroup on $E$ given by the formula
\begin{equation}
\label{ex-11} 
T(t)x=e^{At}x=\sum_{k=0}^{\infty}\frac{t^kA^kx}{k!} 
\end{equation}
and 
the semigroup $\{T(t)\}_{t\ge 0}$ is 
\textit{uniformly continuous},\index{uniformly continuous semigroup} i.e.,
\begin{equation}
\label{d:uc}
\lim_{t\to t_0}\|T(t)-T(t_0)\|=0 \textrm{ \,for $t_0\ge 0$}.
\end{equation}

We will study mainly semigroups of contractions, so we only characterize the generators of such semigroups. 
We recall that a linear operator $T$ on $E$ is a \textit{contraction}\index{contraction} if $\|Tx\|\le \|x\|$ for $x\in E$.
     
\begin{theorem}[Hille--Yosida]\index{Hille--Yosida theorem}
\label{tw-HY}
Let $E$ be a Banach space, and let $A$ be a linear operator on $E$ with~domain
$\mathfrak D(A)$.
The operator $A$ is a generator of a strongly continuous semigroup of contractions if and~only if the following conditions are satisfied:
\begin{enumerate}
\item[\rm(a)]  the set $\mathfrak D(A)$ is dense in~$E$;
\item[\rm(b)]  for each  $\lambda>0$ there exists an operator $(\lambda I-A)^{-1}$ defined everywhere on $E$;
\item[\rm(c)]  $\|(\lambda I-A)^{-1}\|\le \lambda^{-1}$ for each  $\lambda>0$. 
\end{enumerate}
\end{theorem}

The operator $(\lambda I-A)^{-1}$ is called 
the \textit{resolvent}\index{resolvent} of $A$ at~$\lambda$
and it is denoted by $R(\lambda,A)$ or $R_{\lambda}$.
The resolvent is defined for $\lambda$ from the \textit{resolvent set}\index{resolvent!set} $\rho(A)$ which consists of all $\lambda\in\mathbb C$
such that the operator  $\lambda I-A$ is invertible and the inverse operator is a bounded operator defined on the whole space $E$.
The set $\sigma(A)=\mathbb C\setminus \rho(A)$ is called the \textit{spectrum}\index{spectrum} of the operator $A$. 
If the operator $A$ generates a semigroup $\{T(t)\}_{t\ge 0}$,
then $R_{\lambda}$ is called the
\textit{resolvent of the semigroup}
$\{T(t)\}_{t\ge 0}$. 
The resolvent of a strongly continuous semigroup of contractions
is also given by the formula
\begin{equation}
\label{Hf-res}
R_{\lambda}x=\int_0^{\infty} e^{-\lambda t} T(t)x\,\,dt \quad\ \text{for $x\in E$}.
\end{equation}

\begin{remark}
We can also consider nonlinear strongly continuous semigroups.
A family $\{T(t)\}_{t\ge0}$ of nonlinear 
operators on a Banach space is called
 a \textit{nonlinear strongly continuous semigroup}\index{nonlinear strongly continuous semigroup} 
 if it satisfies conditions (a), (b), and  the map $(t,x) \mapsto T(t)x$ is continuous.
\end{remark}

\subsection{Stochastic and substochastic semigroups}
\label{ss:s-s semigroups}
Let  $(X,\Sigma,m)$ be a $\sigma$-finite measure space
and let $L^1=L^1(X,\Sigma,m)$.
A $C_0$-semigroup $\{P(t)\}_{t\ge0}$ of positive contractions on $L^1$ is called
a \textit{substochastic semigroup}.\index{substochastic semigroup}
A $C_0$-semigroup of stochastic operators is called
a \textit{stochastic semigroup}\index{stochastic!semigroup}
or a \textit{Markov semigroup}.\index{Markov!semigroup}
We can also consider a \textit{nonlinear stochastic semigroup}
if $\{P(t)\}_{t\ge0}$ is a nonlinear strongly continuous semigroup
and $P(t)(D)\subset D$ for $t\ge 0$, where $D$ is the set of all densities.

An example of a stochastic semigroup is a semiflow of Frobenius--Perron operators $\{P^t\}_{t\ge 0}$ considered in Sec.~\ref{ss:ergodic}
if we assume that this semiflow satisfies condition (c). 
We will now examine what conditions a linear operator $A$ on a $L^1$ space should satisfy in order to be a generator of a stochastic semigroup.

First we consider the case when $A$ is a bounded operator.
A bounded linear operator $A$  is a generator of
 a positive semigroup on $L^1$ (generally on  a Banach lattice)
if and only if  $A+\|A\|I\ge 0$
 (see e.g. Theorem 1.11, page 255 \cite{Arendt}).
 We recall that a semigroup $\{P(t)\}_{t\ge 0}$ is \textit{positive}\index{positive!semigroup} if each operator $P(t)$ is positive.
Let $\lambda=\|A\|$ and $B=A+\lambda I$,
then $B$ is a positive operator. Thus the generator of the semigroup of positive operators  $\{P(t)\}_{t\ge0}$
on $L^1$ is of the form $A=-\lambda I+B$, where $B$ is a positive and bounded operator. 
If $A$ is a generator of a stochastic semigroup  
then $P(t)$ preserves  the integral.  
Since $P(t)=I+tA+o(t)$, we obtain 
\begin{equation}
\label{int-zero}
\int_X Af(x)\,m(dx)=0\quad \textrm{for $f\in L^1$}.
\end{equation}
Thus $\int_X Bf(x)\,m(dx)=\lambda$ for $f\in D$.
Assume that   $\lambda>0$ and let  $P=B/\lambda$. Then 
$A=-\lambda I+\lambda P$, where  $P$ is a stochastic operator.
Analogously if a linear and bounded operator $A$ is a generator of a substochastic semigroup,  
then $A$ is of the form   
$A=-\lambda I+\lambda P$, where  $P$ is a substochastic operator and $\lambda\ge 0$. The semigroup  $\{P(t)\}_{t\ge0}$ is given by the formula 
\begin{equation}
\label{kang-wz-s}
P(t)f=\sum_{k=0}^{\infty}\frac {(\lambda t)^ke^{-\lambda t}}{k!}P^kf.
\end{equation}
 
\begin{example}
\label{ex:Markov-chain-n} 
Let $X=\{1,2,\dots,n\}$, $\Sigma=2^X$ and let $m$ be the counting measure. In this case the space 
$L^1(X,\Sigma,m)$ is denoted by $l^1_n$.
In the space $l^1_n$ a function $f\colon X\to\mathbb R$ is represented as a sequence $y=(y_1,\dots,y_n)$,  
the integral of $y$ over $X$ is given by
$\sum_{i=1}^n y_i$, and the space
$l^1_n$ is isomorphic to the space  $\mathbb R^n$ with the norm $\|y\|=|y_1|+\dots+|y_n|$. 
The generator $A$ of a stochastic semigroup on $l^1_n$ is of the form
$Ay=yQ$, where the matrix $Q=[q_{ij}]$ satisfies 
conditions 
\begin{enumerate}
\item[(i)] $q_{ij} \ge 0$  for $i\not= j$,
\item[(ii)] $\sum_{j=1}^n q_{ij}=0$ for $i=1,\dots,n$.
\end{enumerate}
The matrix $Q$ is the \textit{intensity matrix}\index{intensity matrix} of a Markov chain $(\xi_t)_{t\ge 0}$.
The entry $q_{ij}$, $i\ne j$, is the rate of jump from the state $i$ to $j$, precisely  
\begin{equation}
\label{inten-matrix}
{\rm P}(\xi_{t+\Delta t}=j|\xi_t=i)=q_{ij}\Delta t+o(\Delta t).
\end{equation} 
Moreover,  $|q_{ii}|$  is the rate of leaving the state $i$.
\end{example}

\begin{remark}
\label{r:Af=0}
One can ask the question whether condition~\eqref{int-zero} is sufficient
 for a bounded linear operator to generate a stochastic semigroup.
 The answer is negative. Let $X=[0,1]$ and $m$ be the Lebesgue measure.     Consider the operator $A=I-P$, where $Pf(x)=2f(2x)$ for $x\in[0,1/2]$ and
 $Pf(x)=0$ otherwise. The operator $A$ satisfies~\eqref{int-zero} and 
generates a semigroup 
$\{T(t)\}_{t\ge0}$ but this semigroup is not positive.
We have
\[
T(t)f=e^te^{-tP}f=e^t\sum_{k=0}^{\infty}\frac {(-1)^kt^kP^k}{k!}f.
\] 
If $f=\mathbf 1_{(1/2,1]}$, then $P^kf=2^k\mathbf 1_{(1/2^{k+1},1/2^k]}$
and $T(t)f(x)<0$ for $x\in (1/2^{k+1},1/2^k]$ and $k$ odd.
\end{remark}

Theorem~\ref{tw-HY} can be applied to check whether an operator $A$ is a generator of a stochastic or substochastic semigroup.
From (\ref{Hf-res}) it follows that if a semigroup is substochastic, then the operators 
$R_{\lambda}$ are positive for each $\lambda>0$. 
It turns out that vice versa, too, if $A$ satisfies conditions (a)--(c) of Theorem~\ref{tw-HY}
and
\begin{enumerate}
\item[(d)] $R_{\lambda}f\ge 0$ for all $\lambda>0$ and $f\ge 0$, 
\end{enumerate}
then $A$ is the generator of a substochastic semigroup.

If a semigroup $\{P(t)\}_{t\ge 0}$ is stochastic, then from  (\ref{Hf-res}) it
follows that the operators $\lambda R_{\lambda}$ are stochastic. If the 
operator $A$ satisfies conditions (a)--(d) and additionally 
\begin{enumerate}
\item[(e)]  for some $\lambda>0$ the operator $\lambda R_{\lambda}$ is stochastic,
\end{enumerate}
then $A$ generates a stochastic semigroup. Condition (e) can be replaced by the following:
\begin{enumerate}
\item[(e${}'$)]   $\int_X Af(x)\,m(dx) =0$ for $f\in \mathfrak D(A)$.
\end{enumerate}
 
In applications, the following problem often arises.
We know that an operator $A$ is the generator of a stochastic or substochastic semigroup $\{P(t)\}_{t\ge 0}$ and $B$ is a linear operator on $L^1$.  
When the operator $A+B$  also generates a stochastic  or a substochastic semigroup? This question is answered by perturbation theorems. 
Now we present one of them useful in applications.  

\begin{theorem}
\label{tw-DP}
Let $A$ be a generator of a stochastic semigroup 
$\{S(t)\}_{t\ge 0}$,
$K$ be a stochastic operator and $\lambda>0$.
The semigroup $\{P(t)\}_{t\ge 0}$ generated by the operator $A+\lambda K-\lambda I$ is stochastic and it  is given by the Dyson--Phillips expansion:
\begin{equation}
\label{eq:dpf1b}
P(t)f=e^{-\lambda t}\sum_{n=0}^\infty \lambda^n S_n(t)f, 
\end{equation}
where
\begin{equation}
\label{eq:dpf2b}
S_0(t)f=S(t)f,\quad S_{n+1}(t)f=\int_0^tS_{n}(t-s)KS(s)f\,ds, \quad
n\ge 0.
\end{equation}
If $\{S(t)\}_{t\ge 0}$ is a substochastic semigroup and $K$ is a substochastic operator then $\{P(t)\}_{t\ge 0}$ is a substochastic semigroup.
\end{theorem}

\subsection{Semigroups related to Markov chains}
\label{ss:Markov-chains}
It can be expected that the construction of the stochastic semigroup presented in Example~\ref{ex:Markov-chain-n}
can be extended to the case where $X$ is a countable set.
However, this case is more difficult. 

Let $l^1=L^1(\mathbb N,2^{\mathbb N},m)$, where $m$ is the counting measure.
Now elements of $l^1$ are infinite sequences $y=(y_0,y_1,\dots)$ and
$\|y\|= \sum_{i=0}^{\infty}|y_i|$.
We assume that an infinity dimensional matrix $Q=[q_{ij}]$
satisfies conditions 
\begin{enumerate}
\item[(i)]  $q_{ij} \ge 0$  for $i\not= j$,
\item[(ii)] $\sum_{j=0}^{\infty} q_{ij}=0$ for $i=0,1,\dots$.
\end{enumerate}

If $\sup_{i\in \mathbb N} |q_{ii}|<\infty$, then 
the operator $Ay=yQ$ is bounded and $P(t)y=e^{At}y$, $t\ge 0$, is a stochastic semigroup on $l^1$. 

The case when $A$ is an unbounded operator is more complex. 
Two linear subspaces can be associated with the matrix $Q$: 
\[
 \mathfrak D_0(Q)= \{ x\in l^1\colon \sum_{i=0}^{\infty} |q_{ii}| |x_i|<\infty \}, 
\]
\[
\begin{aligned}
 \mathfrak D(Q)=\bigg\{ x\in l^1\colon  
&\text{ for each  $j\in\mathbb N$
\ the series  $\sum_{i=0}^{\infty} x_iq_{ij}$ is absolutely} \\ 
 &\text{ convergent and \,}\sum_{j=0}^{\infty}\bigg|\sum_{i=0}^{\infty} x_iq_{ij}\bigg| <\infty \bigg\}. 
\end{aligned}
\]
Firstly, a mapping $A$ with domain $\mathfrak D_0(Q)$ or $\mathfrak D(Q)$
need not be the generator of a strongly continuous semigroup on $l^1$.
Secondly, even when $A$ with an appropriately chosen domain generates a strongly continuous semigroup on $l^1$, 
this semigroup may not be a stochastic semigroup. Thirdly, the same matrix $Q$ 
can be used to generate different positive semigroups (including stochastic semigroups)
on $l^1$ depending on the choice of the domain of $A$.

We consider a family $\mathfrak F$ of all positive semigroups 
on $l^1$ 
with generators $A$ such that: 
$\mathfrak D_0(Q)\subseteq \mathfrak D(A)$ 
and $Ax=xQ$ for $x\in \mathfrak D_0(Q)$.
A semigroup  $\{P(t)\}_{t\ge 0}$ from $\mathfrak F$  is called \textit{minimal}\index{minimal semigroup}
if for any other semigroup $\{T(t)\}_{t\ge 0}$ from $\mathfrak F$  we have
 $T(t)f\ge P(t)f$ for each nonnegative $f\in \mathfrak D_0(Q)$ 
and $t\ge 0$.
\begin{theorem}
\label{c:sub}
Assume that the matrix $Q$  satisfies {\rm (i)--(ii}).
Then there is the minimal  semigroup $\{P(t)\}_{t\ge 0}$ and 
$\{P(t)\}_{t\ge 0}$ is a substochastic semigroup.
\end{theorem}

The answer for the question when a minimal semigroup is stochastic is given by following Kato theorem~\cite{kato54}.
\begin{theorem}
\label{Kol} 
Let  $\lambda>0$ be a positive constant. The minimal semigroup
related to $Q$ is a stochastic semigroup if and only if the equation 
$Qx=\lambda x$ has no nonzero and nonnegative solution $x\in l^{\infty}$. 
\end{theorem}

If the minimal semigroup $\{P(t)\}_{t\ge 0}$ related to 
a matrix~$Q$ is a stochastic semigroup, then
the matrix~$Q$ and the semigroup $\{P(t)\}_{t\ge 0}$ are called \textit
{non-explosive}.\index{non-explosive!matrix}\index{non-explosive!semigroup}
If the matrix is non-explosive, then $Q$ is  
an intensity matrix of a Markov chain $(\xi_t)_{t\ge 0}$ on $\mathbb N$,
i.e. the entries $q_{ij}$ satisfy condition (\ref{inten-matrix}).

\begin{remark}
\label{domain-M-ch}
If $\mathfrak D(A)$ is the domain of the generator of a stochastic semigroup related to the intensity matrix $Q$, then  
\[
\mathfrak D_0(Q)\subseteq \mathfrak D(A)\subseteq \mathfrak D(Q)
\]
and none of these inclusions can be replaced by equality.
\end{remark}

\subsection{Diffusion semigroup} 
\label{ss:d-s}
The theory of stochastic processes provides many interesting examples of
semigroups. Semigroups related to diffusion processes play a fundamental role.

Consider the  \textit{It\^o  equation}\index{It\^o  equation} 
\begin{equation}
\label{IE}
dx_t=\sigma(x_t)\, dw_t + b(x_t)\,dt,
\end{equation}
where $w_t$ is an $m$-dimensional Brownian motion,
$\sigma(x)=[\sigma_j^i(x)]$ is a $d\times m$ matrix and
$b(x)$ is a vector in $\mathbb R^d$  for every $x\in\mathbb R^d$. We assume that
for all $ i=1,..., d$, $j=1,..., m$ the functions
$b _i$, $\sigma_j^i$ are sufficiently smooth and have bounded
derivatives of all orders, and the  functions $\sigma_j^i$
are also bounded. A process  $x_t$, $t\ge 0$, which satisfies (\ref{IE}) is called \textit{diffusion process}.\index{diffusion process}

Assume that the random variable $x_0$  has the density $f(x)$. Then
solution $x_t$ of  (\ref{IE}) 
also  has  the density $u(t,x)$ and $u$ satisfies 
the \textit{Fokker--Planck equation}\index{Fokker--Planck equation}
(also called  the \textit{Kolmogorov forward equation}\index{Kolmogorov forward equation})   
\begin{equation}
\label{F-P}
\frac{\partial u}{\partial t}=
\sum_{i,j=1}^d
\frac{\partial^2(a_{ij}(x)u)}
{\partial x_i\partial x_j}-
\sum_{i=1}^d
\frac{\partial (b_i(x)u)}{\partial x_i},
\end{equation}
where $a_{ij}(x)=\frac{1}{2}\sum_{k=1}^m
\sigma^i_k(x)\sigma^j_k(x)$.
Let $P(t)f(x)=u(t,x)$. Then $\{P(t)\}_{t\ge 0}$ is a stochastic semigroup
(\cite{Rudnicki-LN}).

Consider the operator 
\begin{equation}
\label{gFP}
A_0f=\sum_{i,j=1}^d
\frac{\partial^2(a_{ij}(x)f)}
{\partial x_i\partial x_j}-
\sum_{i=1}^d
\frac{\partial (b_i(x)f)}{\partial x_i}
\end{equation} 
on the set 
$E=\{f\in
L^1(\mathbb R^d)
\cap C^2_b(\mathbb R^d)
\colon 
A_0f\in L^1(\mathbb R^d)\}$, 
where $C^2_b(\mathbb R^d)$ denotes the set of all twice 
differentiable bounded functions whose derivatives of order 
${}\le 2$  are continuous and bounded. 
The closure of the operator $A_0$ is the generator of a stochastic semigroup  
$\{P(t)\}_{t\ge 0}$ (\cite[Proposition 1.3.3]{EK}). 
We recall that 
the \textit{closure  of the operator} $\,A_0\,$\index{closure of a linear operator} is the operator $A$ whose graph is the closure of the graph of $A_0$,
as long as such an operator $A$ exists.

If the functions $a_{ij}$
satisfy the uniform elliptic condition:
\begin{equation}
\label{UE}
\sum_{i,j=1}^d a_{ij} (x)\lambda_i\lambda_j\ge\alpha\|\lambda\|^2
\end{equation}
for some $\alpha>0$ and every $\lambda\in\mathbb R^d$
and $x\in\mathbb R^d$, then
the stochastic semigroup  
$\{P(t)\}_{t\ge 0}$
is an
\textit{integral semigroup},\index{integral!semigroup} i.e.
\[
P(t)f(y)=\int_{\mathbb R^d} q(t,x,y)f(x)\,dx, \quad t>0,
\]
and the \textit{kernel} $q$ is continuous and 
positive.

It should be mentioned that in the case when $\sigma$ is a null matrix,
then Eq.~\eqref{IE} reduces to Eq.~\eqref{M1}. Consequently,   
the generator of the group of Frobenius--Perron operators from Example~\ref{e:F-P-cw} is of the form 
\begin{equation}
\label{gFP-0}
Af=-\sum_{i=1}^d
\frac{\partial (b_i(x)f)}{\partial x_i}.
\end{equation} 

\subsection{Processes with jumps} 
\label{ss:proc-jump}
Let $(X,\Sigma,m)$ be a $\sigma$-finite measure space 
and let $(x_t)_{t\ge 0}$ be a family of Markov processes such that
if the distribution of the random variable $x_0$ has a density $f$,
then  the random variable $x_t$ has a density $S(t)f$ and
$\{S(t)\}_{t\ge 0}$ is a stochastic semigroup. 
Let $(A_0,\mathfrak D(A_0))$ be the generator of the semigroup $\{S(t)\}_{t\ge 0}$.
For example $\{S(t)\}_{t\ge 0}$ can be a diffusion semigroup or a semigroup of Frobenius--Perron operators. 

Consider a family of Markov processes
$(\xi_t)_{t\ge 0}$ 
describing the motion of points in the space $X$
and defined in the following way.
A point moves along the sample paths  of the process
$(x_t)_{t\ge 0}$ but it can randomly jump to a new position according to a jump function defined below and then
it  continues to move along the sample paths of the process $(x_t)_{t\ge 0}$, \textit{etc}.

We define the jump function by specifying the jump intensity from a point $x$ to a set $B\in\Sigma$.
For a fixed point $x\in X$, the 
\textit{jump function}\index{jump!function}
is  a measure $\mathcal K(x,B)$ defined in such a way that the probability of jump from a point $x$ to a set $B$
during the time $\Delta t$ is $\mathcal K(x,B)\Delta t +o(\Delta t)$.
Usually the jump function $\mathcal K(x,B)$ is of the form
$\mathcal K(x,B)=\psi(x)\mathcal P(x,B)$, where $\psi(x)$ is the rate of jump
and $\mathcal P(x,B)$ is the probability of jump from $x$ to a given set $B$. 
The function $\mathcal P$ is 
a \textit{transition probability function}\index{transition probability function} on $(X,\Sigma)$, i.e.  $\mathcal P(x,\cdot)$ 
is a probability measure on $(X,\Sigma)$ for each $x\in X$
and the function
$x\mapsto\mathcal P(\cdot,B)$ is measurable for each $B\in \Sigma$.

Assume that $\mathcal P$ has the following property
\begin{equation}
\label{trana}
m(B)=0 \Longrightarrow 
\mathcal  P(x,B)=0\textrm{ for $\,m$-a.e. $x\,$ and $B\in \Sigma$}.
\end{equation}
Then for every $f\in D$ the measure 
\[
\mu(B)=\int f(x)\mathcal P(x,B)\,m(dx)
\]
is absolutely continuous with respect to the
measure $m$.
This fact is a simple consequence of the
 \textit{Radon--Nikodym theorem},\index{Radon--Nikodym theorem}
 which says that 
 the measure $\nu$ is absolutely continuous 
 with respect to the measure $m$ if and only if
 the following implication $m(B)=0 \Rightarrow \nu(B)=0$ holds for all sets $B\in\Sigma$. 
Now, the formula   
$Pf=d\mu/dm$ defines a stochastic operator $P:L^1\to L^1$ called 
a \textit{jump operator}\index{jump!operator}. 
Moreover, if 
$P^*\colon L^{\infty}\to L^{\infty}$ is the adjoint
operator of $P$, then 
$P^*g(x)=\int g(y)\,\mathcal P(x,dy)$.

We consider two examples of the jump operator. 
One of them is Frobenius--Perron operator $P_{\varphi}$ introduced in Sec.~\ref{ss:ergodic}. The transition probability function corresponding to
$P_{\varphi}$ is of the form
\begin{equation}
\mathcal P(x,B)=
\begin{cases} 
1, \textrm{ if $\varphi(x)\in B$,}\\
0, \textrm{ if $\varphi(x)\notin B$.}
\end{cases}
\end{equation}

The second example is an integral stochastic operator.
If $(X,\Sigma,m)$ is a $\sigma$-finite measure space
 $q\colon X\times X\to[0,\infty)$  is 
a measurable function 
such that
\[
\int_X q(x,y)\,m(dy)=1
\]
for almost all $x\in X$, then
\begin{equation}
\label{integ}
Pf(y)=\int_X q(x,y)f(x)\,m(dx)
\end{equation}
is a stochastic operator. 
The function $q$ is called the {\it kernel} of the operator $P$
and  $P$ is called an
\textit{integral} or \textit{kernel} 
operator.\index{integral!operator}\index{kernel!operator} 
The function $\mathcal P(x,B)=\int_B q(x,y)\,m(dy)$ is
the transition probability function corresponding to $P$.

We assume that the function $\psi\colon X\to [0,\infty)$
describing the intensity of jumps is measurable and bounded.
Then the process $(\xi_t)_{t\ge 0}$ defines a semigroup 
$\{P(t)\}_{t\ge 0}$ on $L^1$ with generator 
\begin{equation}
\label{generator-A0}
Af=A_0f- \psi f+ P(\psi f).
\end{equation}
We check that $\{P(t)\}_{t\ge 0}$ is a stochastic semigroup.
Let $\lambda=\sup_X\psi(x)$  and 
$Kf=\lambda^{-1}[P(\psi f)+(\lambda-\psi)f]$. Then $K$ is a stochastic operator
and $A=A_0-\lambda I+\lambda K$. According to Theorem~\ref{tw-DP} the operator $(A,\mathfrak D(A_0))$ is the generator of a stochastic semigroup
defined by formulas \eqref{eq:dpf1b} and \eqref{eq:dpf2b}.

If $A_0$ is of the form \eqref{gFP-0} and is therefore the generator of a semigroup of Frobenius--Perron operators, then $(\xi_t)_{t\ge 0}$ is an example 
of a piecewise deterministic process~\cite{davis84}. 
A stochastic process $(\xi_t)_{t\ge 0}$ 
is called a \textit{piecewise deterministic processes}\index{piecewise deterministic process} if  there is an increasing sequence of random times $(t_n)$, called jumps, such that sample paths  of $(\xi_t)_{t\ge 0}$
are defined in a deterministic way in each interval $[t_n,t_{n+1})$;
in our case they satisfy Eq.~\eqref{M1}.

\subsection{Diffusion process with random switching}
\label{ss:diff-with-switching}
\begin{figure}
\begin{center}
\begin{picture}(260,200)(0,0)
\includegraphics[scale=0.8,viewport=142 462 456 226]{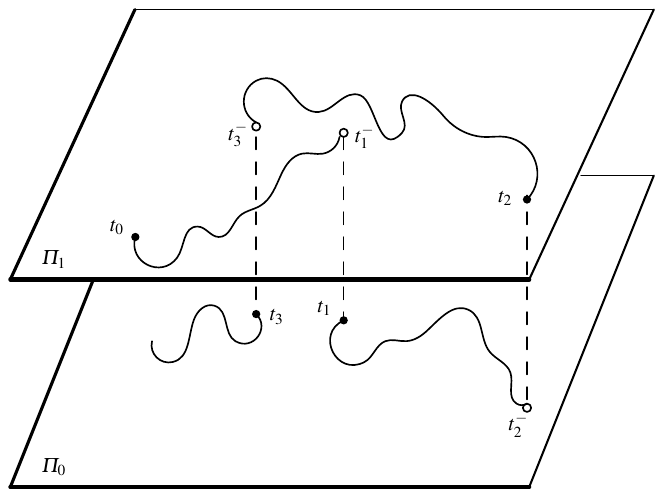}
\end{picture}
\end{center}
\caption{Example of a trajectory of a diffusion process with random switches. In this case $X=\mathbb R^2\times\{0,1\}$ and $t_0,t_1,t_2,t_3$ 
are the successive switching times.}
\label{r:dyf-z-przel}
\end{figure}

Consider the following stochastic equation:
\begin{equation}
\label{IE-s}
dx_t=\sigma(x_t,y_t)\, dw_t + b(x_t,y_t)\,dt,
\end{equation}
where the process $(y_t)_{t\ge 0}$ takes values in the set
$I=\{0,1,\dots,n\}$, $n\ge 1$. We assume that the matrices $\sigma(x,k)$ 
and the vectors $b(x,k)$ satisfy the same assumptions as $\sigma(x)$
 and $b(x)$ from Sec.~\ref{ss:d-s}.  
We also assume that if  $y_t=k$,  then 
the process $(y_t)_{t\ge 0}$ jumps to the state $l$ 
with rate $q_{kl}(x_t)$, i.e. 
\begin{equation}
\label{IE-s1}
{\rm P}(y_{t+\Delta t}=l|y_t=k)=q_{kl}(x_t)\Delta t+o(\Delta t).
\end{equation} 
We assume that the functions $q_{kl}$ are bounded, continuous and nonnegative. 
Then $\xi_t=(x_t,y_t)$, $t\ge 0$ is a Markov process with
values in the space $X=\mathbb R^d\times I$
and it is called 
\textit{diffusion process with random switching} \cite{Y-Zhu}.\index{diffusion process!with random switching} 
Fig.~\ref{r:dyf-z-przel} shows an example of the trajectory of such a process.

If $I=\{0,1\}$ and $\sigma(x,1)=0$, then Eq.~\eqref{IE-s}
describes the process which randomly jumps
between stochastic and deterministic states.
Such processes appear in transport phenomena in sponge--type
structures \cite{BBH,Be,LR}. 

We will now find the stochastic semigroup related to the process  $(\xi_t)_{t\ge 0}$.
Let us denote by $\{S^k(t)\}_{t\ge0}$ the stochastic semigroup with diffusion matrix $\sigma(x,k)$ and with drift vector  $b(x,k)$.
It has the generator of the form
\begin{equation}
\label{gFP-k}
A_kf=\sum_{i,j=1}^d
\frac{\partial^2(a_{ij}(x,k)f)}
{\partial x_i\partial x_j}-
\sum_{i=1}^d
\frac{\partial (b_i(x,k)f)}{\partial x_i}.
\end{equation} 
Let $q_{kk}(x)=-\sum_{l\ne k}q_{kl}(x)$.
For a function $f$ of the variables  $(x,k)$ we use the notation  $f=[f_0,\dots,f_n]$, where $f_k(x)=f(x,k)$.
Let  $Af=[A_0f_0,\dots, A_n f_n]$ and $Qf=[f_0,\dots,f_n]Q$. 
Then the operator $A+Q$ is the generator of a stochastic semigroup 
$\{P(t)\}_{t\ge0}$ associated with the process
$(\xi_t)_{t\ge 0}$.
This semigroup is defined on the space $E=L^1(X,\mathcal B(X),m)$,
where  $m$ is the product of the  Lebesgue measure on $\mathbb R^d$ and the counting measure on  $I$.
The operator $A$ is the generator of a stochastic semigroup 
$\{S(t)\}_{t\ge0}$ on  $E$ given by the formula 
\[
S(t) f=(S^0(t)f_0,\dots, S^n(t)f_n).
\]
Let $\lambda=\sup_{x\in \mathbb R^d} \max\,\{-q_{00}(x),\dots,-q_{nn}(x)\}$ and 
$K=\lambda^{-1}Q+I$. Then  $K$ is a stochastic operator 
and $Q= -\lambda  I+\lambda K$. Since the generator of the semigroup $\{P(t)\}_{t\ge0}$ is of the form  $A-\lambda  I+\lambda K$, 
this semigroup is given by the formulae  
(\ref{eq:dpf1b}) and (\ref{eq:dpf2b}).

In the case when $\sigma\equiv 0$, the Markov process 
$(\xi_t)_{t\ge 0}$ is a piecewise deterministic.
This process describes the action of several \textit{dynamical systems with random switching}\index{dynamical systems with random switching} between them.
Such dynamical systems have many applications in biology. 
They  are  natural models for describing gene regulatory networks.
Depending on the state of activity of the genes present in the network, it is described by different systems of differential equations.
A change in gene activity switches the system~\cite{BLPR,lipniacki2,RT,SH,ZFL,ZFML}.

\subsection{Processes with non-random moments of jumps}
\label{ss:proc-non-random-jumps}
We will now present some processes in which a change in dynamics or a jump in~phase space occurs when the process 
reaches the boundary of the domain  or some fixed subset of the phase space.
An example of such a process is stochastic billiards \cite{CPSV,Evans-b,LM-KR,M-KR}. Similar properties have  biological models of 
cell cycle~\cite{LeR,PR-cell-cycle,PR-cykl-kom-prz,Rotenberg}, neuron activity~\cite{Burkitt,PR-Stein,Stein1967,SGJ}, gene regulatory networks~\cite{HWS,KKLL}, and vesicular transport in the cell~\cite{BN1,BN2,BN3}.

\begin{figure}
\begin{center}
\begin{picture}(360,80)(0,0)
\includegraphics[scale=0.7,viewport=182 562 382 456]{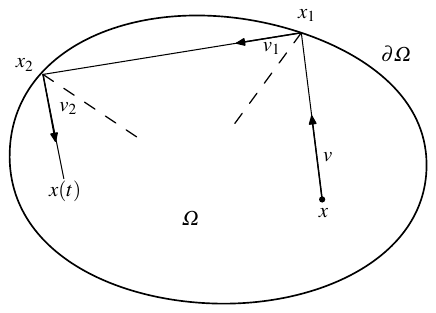}\includegraphics[scale=0.8,viewport=162 592 322 486]{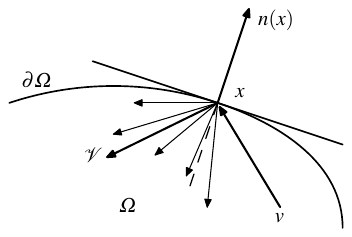}
\end{picture}
\end{center}
\caption{Left: classical billiards -- angles of incidence are equal to the angles
of reflection. Right: stochastic billiards. Regular and diffuse reflection: 
 $\mathcal V$ -- the specular reflection, thin vectors -- diffuse reflection.}
\label{p:col-3}
\end{figure}

Let  $\Omega\subset \mathbb R^d$ be a compact set with nonempty interior and 
with piecewise smooth boundary $\partial \Omega$. 
In \textit{stochastic billiards}\index{stochastic!billiards} a moving point 
does not change its velocity in the interior of $\Omega$,
but when it collides with its boundary $\partial \Omega$. The new velocity is chosen randomly 
from among the velocities directed towards the interior of the set of $\Omega$ and the motion continues until the next collision with the boundary
(see Fig.~\ref{p:col-3}). 

Let $x_t$ be the position of a moving particle  and $v_t$ be its velocity.
Consider a stochastic process $\xi_t=(x_t,v_t)$, $t\ge 0$. This process takes values in the space $X=\Omega\times V$, $V=\mathbb R^d$. 
We can also consider stochastic billiards when a particle is moving with
constant speed $\|v\|=1$. Then $V=S^{d-1}$, where 
$S^{d-1}$ is a unit sphere in the space $\mathbb R^d$.
Assume that if the particle
strikes the boundary at point $x$ and with velocity $v$,
then a new velocity $\bar v$ is chosen according to the density distribution
$k(x,v,\bar v)$.
Assume that the distribution of the random variable $\xi_0$ has density 
$u_0\in L^1(X,\mathcal B(X),\ell)$, where $\ell$ is the Lebesgue measure on $X$.
Then each random variable $\xi_t$ has density $u(t,x,v)$ and the function $u$ satisfies the following initial-boundary problem~\cite{LM-KR}:
\begin{align}
\label{bil-1}
&\frac{\partial u}{\partial t}+v\cdot \nabla_xu(t,x,v)=0,\quad\text{where \ }
\nabla_x=\left[\dfrac{\partial}{\partial x_1},\dots,\dfrac{\partial}{\partial x_d}\right],\\
\label{bil-2}
&u(0,x,v)=u_0(x,v),\\
\label{bil-3}
&u_{|\Gamma_-}=H(u_{|\Gamma_+}),
\end{align}
where 
\[
\Gamma _{\pm }=\left\{ (x,v)\in \partial \Omega \times V;\ \pm
v \cdot n(x)>0\right\},
\]
$n(x)$ is the outward unit normal at $x\in \partial \Omega$
and $H$  is
a linear boundary operator relating the
outgoing and incoming fluxes $u _{\mid \Gamma _{+}}$ and $u_{\mid
\Gamma _{-}}$.
The operator $H$ 
 is linear and bounded  on
the trace spaces
\[
L^{1}_{\pm}=L^{1}(\Gamma_{\pm}\,,\,
|v\cdot n(x)| \pi(dx)\otimes m_V(dv))=L^{1}(\Gamma_{\pm},d\mu_{\pm}(x,v)),
\]
where $\pi$ denotes the Lebesgue surface measure on $\partial \Omega$ 
and $m_V$ is the Lebesgue measure on the space  $V$.
The operator $H\colon L^1_+\to L^1_-$ is given by the formula
\[
H\varphi(x,\bar v)=\int_{\Gamma _{+}(x)}k(x,v,\bar v)
\varphi (x,v)|v\cdot n(x)|\,m_V(dv), 
\]
where $\Gamma _{+}(x)=\{v\in V\colon (x,v)\in \Gamma _{+}\}$.
Under some technical assumptions the solutions of the system 
\eqref{bil-1}--\eqref{bil-3} generate a stochastic semigroup 
on the space $L^1(X,\mathcal B(X),\ell)$
and this semigroup has an interesting long-time behaviour (see~\cite{LM-KR}).

One-dimensional version of stochastic billiards was studied in~\cite{M-KR}.
This version was applied in  \cite{PR-cykl-kom-prz} to study Lebowitz and Rubinow model\index{Lebowitz and Rubinow model} of the cell cycle~\cite{LeR}.
Some models of the cell cycle will be presented in more detail in Sec.~\ref{ss:cykl-kom}.

\section{Long-time behaviour of  semigroups of operators}
\label{s:long-time-beh}
\subsection{The Perron--Frobenius theorem}
\label{ss:gen-case}
A good starting point to study long-time behaviour of semigroups of operators is 
a continuous version of the classical Perron--Frobenius theorem.\index{Perron--Frobenius theorem} 
To recall it, let  $E$ be the space $\mathbb R^n$ with the scalar product $\langle \cdot,\cdot\rangle$. 
The zero vector in $E$ is denoted by $\mathbf 0$.
If $x\in E$, then we write $x>\mathbf 0$ when
$x_i>0$ and $x\ge \mathbf 0$ when $x_i\ge 0$  for $i=1,2,\dots,n$. Let $Q=[q_{ij}]$ be an $n\times n$ real matrix. 
Then the formula $T(t)x=xe^{tQ}$ defines a semigroup $\{T(t)\}_{t\ge 0}$ on $E$.
We assume that 
\begin{enumerate}
\item [(P)]  $q_{ij}\ge 0$ for $i\ne j$,  
\item [(I)]  for all $i$ and $j$  there exists a sequence
$(i_1,i_2,\dots,i_m)$ such that  
$i_1=i$, $i_m=j$ and 
$q_{i_{k}i_{k+1}}>0\,$ for $k=1,\dots,m-1$.
\end{enumerate} 
Condition (P) guarantees that the semigroup $\{T(t)\}_{t\ge 0}$ is positive,
 i.e. $T(t)x\ge \mathbf 0$ for $x\ge \mathbf 0$. 
From both conditions (P) and (I) it follows that there are 
a constant $r\in\mathbb R$  and vectors $x^*>\mathbf 0$ and $y^*>\mathbf 0$ such that
\begin{equation}
\label{c-d2}
\lim_{t\to\infty} e^{-r\, t}T(t)x= x^*\,\langle  y^*,  x \rangle \quad\textrm{for $x\in E$} 
\end{equation}
and the convergence in \eqref{c-d2} is exponential. 
The \textit{irreducibility}\index{irreducibility} assumption (I) implies that
the limit set in~\eqref{c-d2} is one-dimensional.
A semigroup $\{T(t)\}_{t\ge 0}$ satisfying condition \eqref{c-d2} 
is said to have \textit{asynchronous exponential growth}.\index{asynchronous exponential growth}
It should be mentioned that one can find in literature, e.g. \cite{Webb-a}, a more general definition of asynchronous exponential growth,
where it is  merely assumed that $e^{-r\,t}T(t)x$ converges to a nonzero finite rank operator.

This theorem is a consequence of quite nontrivial properties of the matrix $Q$.
First of all, the number $r$ is a \textit{dominant eigenvalue}\index{dominant eigenvalue} of $Q$, which means that any other eigenvalue $\lambda$ satisfies 
$\mathrm{Re\,} \lambda<r$.  Also, the general eigenspace associated to $r$ is one-dimensional, and 
 $x^*$ and $y^*$ are left and right eigenvectors corresponding to $r$, i.e. $x^*Q=rx^*$ and $Qy^*=ry^*$ (it follows that $T(t)x^*=e^{r\,t}x^*$). Finally, the functional $\alpha$ defined by $\alpha(x)= \langle  y^*,  x \rangle$ is positive, i.e.  $\alpha(x)\ge 0$ for $x\ge 0$.

If assumption (I) is dropped, the dominant eigenvalue $r\in \mathbb R$ still exists, but \eqref{c-d2} does not hold. More specifically, if the 
general eigenspace associated to $r$ is one-dimensional  then (P) implies $\eqref{c-d2}$  but
the inequalities $x^*>\mathbf 0$ and $y^*>\mathbf 0$ need not be fulfilled. If  the general eigenspace associated to $r$ is not one-dimensional the behaviour of 
the semigroup $\{T(t)\}_{t\ge 0}$ is more complex.
For example, if $Q$ is the $2\times 2$ matrix with $q_{11}=q_{12}=q_{22}=1$ and $q_{21}=0$, then $r=1$ is the dominant eigenvalue and  
$T(t)x =e^{t}(x_1,x_2+x_1t)$. This  means that $\lim_{t\to\infty} e^{-t}t^{-1}T(t)x=xR$, where $xR=(0,x_1)$.
Generally, if (P) holds, then 
the semigroup $\{T(t)\}_{t\ge0}$ has polynomial exponential growth, i.e.
there exists a nonnegative integer $k$ and a  $n\times n$  nonzero matrix $R$ with nonnegative elements such that, for each $x$, 
$\lim_{t\to\infty} e^{-r\,t}t^{1-k}T(t)x= xR$. More precisely, $k\times k$ is the size of the largest  Jordan block related to the eigenvalue $r$, and $xR$ is a linear combination of certain eigenvectors related to these Jordan blocks.

If (P) is not satisfied, a dominant eigenvalue, defined as an eigenvalue with the largest real part, may be complex, and there may exist a number of eigenvalues with this property. These two new features result in two new types of asymptotic behaviour: using  the Jordan--Chevalley decomposition of $Q$ we show that (a) if imaginary parts of all dominant eigenvalues are commensurable, then the semigroup has a periodic polynomial exponential growth and (b) if they are not, the semigroup has an almost periodic exponential growth. 
To summarise: asymptotic behaviour of the semigroup $\{T(t)\}_{t\ge 0}$ is determined by its behaviour on the invariant subspace 
$E_d$, spanned by generalised eigenspaces, with the largest dimension, related to dominant eigenvalues.

If we turn to semigroups in infinite-dimensional spaces, new problems  
appear because spectral properties of infinite-dimensional semigroups are more complex.
For example, consider the intensity matrix $Q$ indexed by integers, with elements
$q_{i\,i+1}=q_{i\,i-1}=1$, $q_{ii}=-2$ and $q_{ij}=0$ in other cases. Since this $Q$ is related to a bounded linear operator on the space $l^1$, we may think of the semigroup $\{T(t)\}_{t\ge 0}$ defined by $T(t)=e^{tQ}$. It can be shown that $\{T(t)\}_{t\ge 0}$ is a stochastic semigroup, i.e.
$\|T(t)x\|=\|x\|$ for $x\ge 0$. The matrix $Q$ also satisfies conditions (P) and (I) and thus one could expect that there exists $x^*> \mathbf 0$ such that 
$T(t)x^*=x^*$, but that is not the case. In other words, for any $t$, $r=1$ is not an eigenvalue of $T(t)$, even though the spectral radius of $T(t)$ is $1$, and thus 
 condition \eqref{c-d2} fails to hold. The main reasons for this is non-compactness of trajectories $t\mapsto T(t)x$.

The second problem in infinite-dimensional spaces is related to the rate of convergence in \eqref{c-d2}. Let  $\{T(t)\}_{t\ge 0}$ be a semigroup on a Banach space $E$.
We say that this semigroup has \textit{asynchronous} (or \textit{balanced}) \textit{exponential growth}\index{asynchronous exponential growth}
 if there exist $r\in\mathbb C$,  nonzero $x^*\in E$, and a nonzero linear functional $\alpha\colon E\to \mathbb C$ 
such that 
\begin{equation}
\label{c-d3}
\lim_{t\to\infty}e^{-r\,t}T(t)x=x^*\alpha(x)\quad\textrm{for $x\in E$}.
\end{equation}

This definition implies that $r$ is a simple eigenvalue of the infinitesimal generator $A$ of the semigroup $\{T(t)\}_{t\ge 0}$. 
In finite-dimensional cases, convergence in \eqref{c-d3} is exponential because the spectrum of $A$ is finite and all
eigenvalues $\lambda$ of $A$,  other than $r$, satisfy the inequality   $\mathrm{Re\,} \lambda<\mathrm{Re\,} r$.
But in the infinite-dimensional case,  condition \eqref{c-d3} does not guarantee that $r$ is an isolated element of the spectrum of $A$
and the rate of convergence in  \eqref{c-d3} needs not be exponential.

Finally, note that  if $ E_0=\mathrm{Ker\,} \alpha$, then  $ E_0$ is invariant with respect to the semigroup. By restricting $\{T(t)\}_{t\ge 0}$ to this subspace,
the problem of the rate of convergence \eqref{c-d3} is replaced by that of the rate of convergence of 
$\| e^{-r\,t}T(t)x\|\to 0\,$ for $x\in  E_0$.

\subsection{Asymptotic decomposition via spectral properties}
\label{ss:expon-grow}
There are several results on asymptotic behaviour of strongly continuous semigroups which are in the spirit of the previous section. 
These results are based on compactness properties of semigroups and are thoroughly discussed  in monographs \cite{EngelNagel,Nagel}.

Let $A$ be the infinitesimal generator of a strongly continuous semigroup $\{T(t)\}_{t\ge 0}$ on a Banach space $ E$
and let $\mu \in \mathbb C$ be an isolated point of  the spectrum of $A$. Since the resolvent $R_{\lambda}$ is a holomorphic function of $\lambda$,
it can be expanded as a Laurent series
\[
R_{\lambda} =\sum_{n=-\infty}^{\infty}(\lambda-\mu)^nU_n
\]	
for $0 < |\lambda - \mu| < \delta$ and some sufficiently small $\delta > 0$. The coefficients $U_n$ of this
series are bounded operators given by the formulas
\[
U_n =
\frac 1{2\pi i}\oint_{\gamma}\frac{R_{\lambda}}{(\lambda-\mu)^{n+1}}\,d\lambda,\quad n\in\mathbb Z,
\]
where $\gamma$ is the positively oriented boundary of a disc with radius $<\delta$ and center at $\mu$. The operator $P=U_{-1}$ is  called the 
\textit{residue}\index{residue of resolvent} 
of $R_{\lambda}$ at $\mu$. We have $U_{-(n+1)}= (A - \mu I)^nP$.
If there exists a $k > 0$ such that $U_{-k} \ne 0$ while $U_{-n} = 0$ for all $n > k$,
then $\mu$ is called a pole of $R_{\lambda}$ of order $k$.

Let $ (L( E),\|\cdot\|)$ be the space of all bounded linear operators on $ E$ 
endowed with the operator norm. Denote by $K( E)\subseteq L( E)$ the set of all compact operators on $ E$. 
A strongly continuous semigroup $\{T(t)\}_{t\ge 0}$ on   $ E$ is called 
\textit{eventually compact}\index{eventually compact semigroup}
 if there exists $t_0 > 0$ such that the operator  $T(t_0)$ is compact.
It can be shown that an eventually compact semigroup is norm continuous for $t\ge t_0$.
A strongly continuous semigroup $\{T(t)\}_{t\ge 0}$  on   $ E$ is called
\textit{quasi-compact}\index{quasi-compact semigroup} if
\[
\lim_{t\to\infty} \inf\{\|T(t)-K\|\colon K\in  K( E)\}=0.
\]
It is clear that an eventually compact semigroup is quasi-compact but not vice versa.
The asymptotic behaviour of a quasi-compact  semigroup is fully characterized by the following theorem (see \cite[p. 333]{EngelNagel}).

\begin{theorem}
\label{th-as-dec-sp}
Let  $\{T(t)\}_{t\ge 0}$  be a quasi-compact  semigroup  on a Banach space $ E$
with generator $A$. Then 
the set $\{\lambda\in \sigma(A)\colon \mathrm{Re\,}\lambda\ge 0\}$ is finite (or empty) and consists of
poles of $R_{\lambda}$ of finite algebraic multiplicity.
If we denote these poles by $\lambda_1,\dots,\lambda_m$, and their residues by $P_1,\dots,P_m$ with
poles of order $k_1,\dots,k_m$, we have  $T(t) = T_1(t) + T_2(t) +\dots+ T_m(t) + R(t)$, where
\begin{equation}
\label{decomp-v-sp1}
T_n(t) = e^{\lambda_nt}
\sum_{j=0}^{k_n-1}
\frac{t^j}{j!}(A-\lambda_n I)^jP_n,
\quad\textrm{for $t\ge 0$ and $1\le n \le m$},
\end{equation}
and
\begin{equation}
\label{decomp-v-sp2}
\|R(t)\|\le Me^{-\varepsilon t}
 \quad \textrm{for some $\varepsilon>0$, $M\ge 1$ and all $t \ge 0$}.
\end{equation}
\end{theorem}

For a particular instance of this theorem,  consider the case when there exists a unique dominant eigenvalue $r\in \mathbb C$ being a first-order pole. 
Then $r$ is the \textit{spectral bound}\index{spectral bound} $s(A)$ of $A$, i.e.
\[
s(A)= \sup\{\textrm{Re\,}\lambda\colon \lambda\in\sigma(A)\}
\]
and $P=\lim\limits_{\lambda\to r}(\lambda-r)R_{\lambda}$. From Theorem~\ref{th-as-dec-sp} 
it follows that 
\begin{equation}
\label{decomp-v-sp3}
\|e^{-r\,t}T(t)-P\|\le Me^{-\varepsilon t}
 \quad \textrm{for some $\varepsilon>0$, $M\ge 1$ and all $t \ge 0$}.
\end{equation}

Theorem \ref{th-as-dec-sp} is also a proper background for a generalisation of the Perron-Frobenius theorem to positive semigroups on Banach lattices.  
To obtain a version of the Perron--Frobenius theorem for such a semigroup, one typically checks, instead of (I),
that the semigroup is \emph{irreducible};\index{irreducible semigroup} this means that for every 
$x\in  E_+\setminus \{\mathbf 0\}$ and for every $x^*\in E^*_+\setminus \{\mathbf 0\}$ there exists a $t>0$
such that $x^*(T(t)x)>0$ (here $ E_+=\{x\in  E\colon x\ge 0\}$ and $ E^*_+=\{x\in  E^*\colon x\ge 0\}$). Indeed, the irreducibility 
rules out the possibility for the range of $P$ to be more than one-dimensional, and together with
quasi-compactness implies asynchronous exponential growth of the semigroup. 
A similar result can be obtained if we replace quasi-compactness by an assumption based on the notion of the \textit{measure of noncompactness} of $A$, see e.g.  \cite{Webb-a}.

The asynchronous exponential growth \eqref{c-d3} is a special case of the situation when there exists the limit 
\begin{equation}
\label{decomp-v-sp4}
\lim_{t\to\infty} e^{-r\,t}T(t)=P,
\end{equation}
where  $P\colon  E\to  E$ is an arbitrary nonzero operator, and
the convergence can be uniform, strong or weak. The latter property was studied in \cite{Thieme-a,Thieme-b}, 
and  a comprehensive overview of applications of this theory to population dynamics can be found in \cite{Webb-b}.
If $\widetilde T(t)= e^{-r\,t}T(t)$ and $\widetilde{ E}= \mathrm{Ker}\,P$, 
then condition \eqref{decomp-v-sp4} in the case of strong convergence is equivalent to 
\textit{strong asymptotic stability}\index{strong asymptotic stability of a semigroup}
 of the semigroup $\{\widetilde T(t)\}_{t\ge 0}$, i.e. to the condition 
\begin{equation}
\label{decomp-v-sp5}
\lim_{t\to\infty}\|\widetilde T(t)x\|= 0\quad \textrm{for $x\in \widetilde{ E}$}.
\end{equation}
In the case when the semigroup is bounded and the set $\sigma(A)\cap(i\mathbb R)$ is at most countable,
the strong asymptotic stability holds 
if and only if the adjoint operator $A^*$ has no pure imaginary eigenvalues \cite{ArendtBatty,Lyubich-Phong}.
Some versions of this result for unbounded operators are proved in \cite{SklyarPolak}.  

The question of deriving the optimal rate of decay in \eqref{decomp-v-sp5} is highly nontrivial.
This problem is extensively discussed in \cite{BattyChillTomilov}.
To give a taste of this theory we recall one of the results of \cite{BorichevTomilov} (Theorem 2.4). 
Let $\{T(t)\}_{t\ge 0}$ be a strongly continuous semigroup on a Hilbert space $ E$ with generator A such that 
$\sigma(A)\cap(i\mathbb R)=\emptyset$. Then for a fixed $\alpha > 0$ the following conditions are equivalent:
\begin{enumerate}
\item[(i)] $R_{is} = O(|s|^{\alpha})$,\quad $s\to\infty$,
\vskip1mm
\item[(ii)] $\|T(t)A^{-1}\|= O(t^{-\alpha})$, \quad $t\to\infty$,
\vskip1mm
\item[(iii)] $\|T(t)A^{-1}x\|= o(t^{-\alpha})$, \quad $t\to\infty$, $x\in E$.
\end{enumerate}

\subsection{Asymptotics  of stochastic semigroups: introductory results}
\label{ss:stoch-semig-s}
We will consider stochastic semigroups on the space
$L^1=L^1(X,\Sigma,m)$ with a $\sigma$-finite measure $m$. 
The iterates of a stochastic  operator
form a \emph{discrete-time} semigroup,\index{discrete-time semigroup} and by writing $P(t)=P^t$ for these iterates,
we may formulate most of definitions and results for both types of semigroups in a unified manner. 

We will begin by studying the asymptotics of uniformly continuous stochastic semigroups. We recall that such a semigroup is given by the formula
\begin{equation}
\label{kang-wz-s2}
P(t)f=\sum_{k=0}^{\infty}\frac {(\lambda t)^ke^{-\lambda t}}{k!}P^kf,
\end{equation}
where $P$ is a stochastic operator and $\lambda\ge 0$. 
\begin{theorem}
\label{th:zbieznosc}
If $Pf^*=f^*$ for some density $f^*>0$ a.e., then 
for any density $f$ there exists an 
invariant density $\bar f$ such that $\bar f=\lim_{t\to\infty} P(t)f$ in $L^1$.
If additionally, $f^*$ is a unique invariant density of $P$, then $\lim_{t\to\infty} P(t)f=f^*$.
\end{theorem}

The proof of this theorem is given in~\cite{PR-prz-apss-apdMp}. 
In the proof we define a new operator 
$\tilde Pf=P(f\!f^*)/f^*$ on the probability space 
$L^1(X,\Sigma,\mu)$ with $\mu(dx)=f^*(x)\,m(dx)$. Then
$\tilde P$ is a \textit{double stochastic operator},\index{double  stochastic operator} i.e. $\tilde P$ is a stochastic operator and $\tilde P\mathbf 1_X=\mathbf 1_X$.  Next we show that 
Theorem~\ref{th:zbieznosc} 
results from the following lemma
\begin{lemma}
\label{lem:Toeplitz2}
If $P$ is a double stochastic operator, then 
$P\E(f|\Sigma_{\mathrm{inv}})=\E(f|\Sigma_{\mathrm{inv}})$
and
\[
\lim_{t\to\infty}\|e^{t( P-I)}f- \E(f|\Sigma_{\mathrm{inv}})\|=0\quad\textrm{for $f\in L^1(X,\Sigma,\mu)$.} 
\]
The $\sigma$-algebra $\Sigma_{\mathrm{inv}}$ is given by
$\Sigma_{\mathrm{inv}}=\{A\in \Sigma \colon P^*\mathbf 1_A=\mathbf 1_A\}$.
\end{lemma}
 
A stochastic semigroup $\{P(t)\}_{t\ge 0}$ is called
{\it asymptotically stable}\index{asymptotic stability} if  there exists a density $f^*$   such that
\begin{equation}
\label{d:as}
\lim _{t\to\infty}\|P(t)f-f^*\|=0 \quad \textrm{for}\quad f\in D.
\end{equation}
We recall that $D$ denotes the set of all densities.
From \eqref{d:as} it follows that $f^*$ is an \textit{invariant density},\index{invariant!density}
i.e. $P(t)f^*=f^*$ for $t\ge 0$.
If the semigroup $\{P(t)\}_{t\ge 0}$ is generated by an
evolution equation $u'(t)=Au(t)$, then its  
asymptotic stability means that the stationary solution $u(t)=f^*$ is  globally asymptotically stable in the sense of Lyapunov\index{Lyapunov} on the set of densities $D$.
Observe that if there is an invariant density $f^*$, then the stochastic semigroup $\{P(t)\}_{ t\ge 0}$
is asymptotically stable if and only if $\lim_{t\to\infty} \|P(t)f\|=0$ for 
$f\in L^1_0=\{f\in L^1\colon \int_X f(x)\,m(dx)=0\}$.

There are several results concerning asymptotic stability of stochastic semigroups. 
We start with  the \textit{lower function theorem}\index{lower function theorem}
of Lasota and Yorke~\cite{LY}.\index{Lasota}\index{Yorke}
Let $f^-(x)=\max(0,-f(x))$. A function $h\in L^1$,  such that $h\ge 0$ and $h\ne 0$, is called 
a \textit{lower function} for a stochastic semigroup $\{P(t)\}_{t\ge 0}$  if 
\begin{equation}
\label{lfu}
\lim_{t\to\infty}\|(P(t)f-h)^-\|=0\quad\textrm{for every $f\in D$.}
\end{equation}

\begin{theorem}
\label{L-Y}
A continuous or discrete-time stochastic semigroup  is asymptotically stable if and only if it has  a lower function $h$. 
\end{theorem}

This theorem was generalized in many directions. The interested reader is referred to the book~\cite{Emelyanov}, where one can find versions   
formulated in the language of von Neumann\index{von Neumann} algebras and to the paper~\cite{Gerlach-Gluck-ETDS} which provides 
a comprehensive overview of earlier results and presents new theorems for semigroups on Banach lattices including 
a version with lower functions depending on $f$.

\subsection{Partially integral semigroups}
\label{ss:partially-int}
In next two sections we present two theorems  on asymptotic behaviour of 
 stochastic semigroups, that stem from the theory of Harris operators\index{Harris} \cite{Fo,Meyn-Tweedie,Num}.

A stochastic semigroup $\{P(t)\}_{t\ge 0}$
is called {\it partially integral}\index{partially integral semigroup}  if
there exist a $t>0$ and a  measurable function $q(t,\cdot,\cdot)\colon X\times X\to [0,\infty)$, called the 
\textit{kernel of the semigroup},\index{kernel!of the semigroup}
 such that
\begin{equation*}
\int_X\int_X  q(t,x,y)\,m(dx)\,m(dy)>0
\end{equation*} 
and
\begin{equation}
\label{w-eta0}
P(t)f(y)\ge \int_X q(t,x,y)f(x)\, m(dx)\quad \textrm{for $f\in D$.}
\end{equation}

\begin{theorem}
\label{asym-th}
Let $\{P(t)\}_{t\ge 0}$
be a partially integral stochastic semigroup.
We assume that this semigroup 
has invariant density $f^*$ and~there are no other periodic points
in the set of densities.
If $f_*>0$ a.e.,
then the semigroup $\{P(t)\}_{t\ge 0}$
is asymptotically stable.   
\end{theorem}

The proof of Theorem~\ref{asym-th} is given in~\cite{R-b95}.  
If $\{P(t)\}_{t\ge 0}$ is a continuous-time stochastic semigroup, 
then Theorem~\ref{asym-th} can be strengthened significantly~\cite{PR-jmaa2}:

\begin{theorem}
\label{asym-th2} 
Suppose that  a continuous-time partially integral stochastic  semigroup $\{P(t)\}_{t\ge 0}$ has
a unique invariant density $f^*$. If $f^*>0$ a.e., then the semigroup
$\{P(t)\}_{t\ge 0}$ is asymptotically stable.
\end{theorem}
\begin{corollary}
\label{asym-irreducible} 
If a continuous-time partially integral stochastic
semigroup is irreducible and has an invariant density, then it is  asymptotically stable.
\end{corollary}
In order to formulate other corollaries to Theorem \ref{asym-th}, we need two  auxiliary definitions.

We say that a stochastic 
semigroup $\{P(t)\}_{t\ge 0}$
{\it spreads supports\,}\index{spreading supports} if 
for every set $A\in\Sigma$ and
for every $f\in D$ we have
\[
\lim_{t\to\infty} m(\mathrm{supp}\, P(t)f\cap A)=m(A)
\]
and 
\textit{overlaps supports,}\index{overlaping supports}  if 
for every $f,g\in D\,$ there exists $t>0$ such that 
\[
m(\mathrm{supp}\, P(t)f\cap 
\mathrm{supp}\, P(t)g)>0.
\]
\begin{corollary}
\label{colo}
A partially integral stochastic semigroup which spreads
supports and has an invariant density 
is asymptotically stable. 
\end{corollary}
\begin{corollary}
\label{colo-n}
A partially integral stochastic semigroup which overlaps
supports and has an invariant density $f_*>0$ a.e. 
is asymptotically stable. 
\end{corollary}

We will now introduce the notion of sweeping, which, together with asymptotic stability, 
will allow us to give a characterisation of the behaviour of 
a rather wide class of substochastic semigroups.

A stochastic semigroup $\{P(t)\}_{t\ge 0}$ is
called \textit{sweeping}\index{sweeping}
with respect to a set $B\in\Sigma$ if for every
 $f\in D$
\begin{equation*}
\lim_{t\to\infty}\int_B P(t)f(x)\,m(dx)=0.
\end{equation*}
The notion of sweeping was introduced in~\cite{KT} but it is also known as 
\textit{zero type property}.\index{zero type property}
In~\cite{KT} the following theorem was proved. 
\begin{theorem}
\label{KT}
Let $\{P(t)\}_{t\ge 0}$
be an integral stochastic semigroup.
Assume that this semigroup 
has no invariant density,
while for a certain set 
$A\in\Sigma$ the following condition holds: 
\vskip2mm
\noindent {\rm (KT)} there exists a measurable function
$f^*$ such that 
$0<f^*<\infty$ a.e. and $P(t)f^*\le f^*$ for $t\ge 0$,
$f^*\notin L^1$ and $\,\int_A f^*(x) \,m(dx)<\infty$.
\vskip2mm
\noindent Then the
semigroup
$\{P(t)\}_{t\ge 0}$
is sweeping with respect to~$A$.
\end{theorem}

In paper \cite{R-b95} it was shown that Theorem \ref{KT} holds for a
broader class of semigroups than integrals, for example to
overlapping semigroups.

\begin{theorem}
\label{KT1}
Let $\{P(t)\}_{t\ge 0}$
be a stochastic semigroup which overlaps supports. 
Assume that 
the semigroup 
$\{P(t)\}_{t\ge 0}$ and a set 
$A\in\Sigma$ satisfy condition {\rm(KT)}.
Then the
semigroup
$\{P(t)\}_{t\ge 0}$
is sweeping with respect to~$A$.
\end{theorem}

The advantage of the theorems~\ref{KT} and~\ref{KT1} is that 
knowing the function $f^*$ we can determine the sets from which the semigroup is sweeping. The main problem 
in applying the theorems~\ref{KT} and~\ref{KT1} 
is to prove that a stochastic semigroup satisfies condition (KT).

\subsection{Asymptotic decomposition via kernel minorants}
\label{ss:r-as-decomp}
Theorems~\ref{asym-th} and~\ref{asym-th2}  allow us to prove the asymptotic stability of Markov semigroups when we know that there is a positive 
invariant density. Unfortunately, in most applications it is difficult to determine the invariant density or even to prove its existence.

Our immediate goal will be to present results on the asymptotic behaviour of substochastic semigroups, 
in which we do not directly assume that the semigroup has invariant density.
One of the corollaries will be a criterion of sweeping from compact sets, not based
 on condition (KT) and convenient for applications.

We present  a method  developed in \cite{PR-JMMA2016,PR-SD}  to study  asymptotic properties of stochastic semigroups via their asymptotic decomposition.
From now on we assume additionally  that $(X,\rho)$ is a separable metric
space and   $\Sigma=\mathcal B(X)$ is the $\sigma$-algebra of Borel subsets of $X$, and consider a partially integral substochastic semigroups $\{P(t)\}_{t\ge 0}$ 
with a kernel $q$ which satisfies the following condition:

\begin{enumerate}
\item[(K)] for every $x_0\in X$ there exist  an $\varepsilon >0$,  a $t>0$,
and a measurable function 
$\eta\ge 0$ such that $\int \eta(x)\, m(dx)>0$ and
\begin{equation}
\label{w-eta3}
q(t,x,y)\ge \eta(y) \quad \textrm{for $x\in B(x_0,\varepsilon)$, $y\in X$},
\end{equation}
where
$B(x_0,\varepsilon)=\{x\in X:\,\,\rho(x,x_0)<\varepsilon\}$.
\end{enumerate}
It is clear that substochastic semigroups satisfying condition (K) are partially integral. 

Now we present two theorems on asymptotic decomposition: the first for a continuous-time substochastic semigroups and the second for substochastic operators 
(discrete-time semigroups) \cite{PR-JMMA2016,PR-SD}. 

\begin{theorem}
 \label{th:2}
Let $\{P(t)\}_{t\ge 0}$ be a continuous-time substochastic semigroup which satisfies {\rm (K)}.
Then there exist an at most countable set $J$, a family of invariant densities
$\{f^*_j\}_{j\in J}$ with disjoint supports $\{A_j\}_{j\in J}$,  and a family
$\{\alpha_j\}_{j\in J}$
of positive linear functionals  defined on $L^1$
 such that
\begin{enumerate}
 \item[\rm (i)] for every $j\in J$ and for every  $f\in L^1$ we have
\begin{equation}
\label{pomoc:ogolny}
 \lim_{t\to \infty}\|\mathbf 1_{A_j} P(t)f-\alpha_j(f)f^*_j\|=0,
\end{equation}
\item[\rm (ii)] if $\,Y=X\setminus \bigcup\limits _{j\in J} A_j$, then for every  $f\in L^1$ and for every compact set $F$ we have
\begin{equation}
\label{th2:sw}
\lim_{t\to \infty} \int\limits_{F\cap Y} P(t)f(x)\, m(dx)=0.
 \end{equation}
 \end{enumerate}
\end{theorem}
We note that not only the sets $A_j$, $j\in J$, are disjoint but also their closures are disjoint, which  
means that the measures $\mu_j(dx)=f^*_j(x)\,m(dx)$ have disjoint topological supports \cite{PR-SD}.

We will precede the formulation of the asymptotic decomposition for substochastic operators with some auxiliary definitions.  

Let $E$ be a vector space. A linear mapping 
 $T\colon E\to E$ satisfying the condition $T^2 = T$ is called a 
\textit{ projection}\index{projection}.
If $y\in T(E)$, then $y=Tx$ for some $x\in E$ and from the definition of projection we get $Ty=T^2x=Tx=y$, so $T$ is the identity on $T(E)$.

Let $(X,\Sigma,m)$ be a $\sigma$-finite measure space. 
A substochastic operator  $S$ is called \textit{periodic} on the set  $B\in\Sigma$, 
if there is a sequence of densities  
$h_1,\dots,h_k$ satisfying the following conditions
\begin{gather}
\label{per-op0}
h_i(x)=0\quad\text{for $i=1,\dots,k$ and $x\notin B$},
\\
\label{per-op1}
h_ih_j=0\quad\text{for $i\ne j$,\quad\  $h_1(x)+\cdots+h_k(x)>0$\quad for $x\in B$,} 
\\
\label{per-op2}
Sh_i=h_{i+1}\quad \text{for $i\le k$, \quad where\quad  $h_{k+1}:=h_1$,} 
\end{gather} 
and  $Sf= ST\!f$ for $f\in L^1(X,\Sigma,m)$, where $T$ is a positive projection
from the space $L^1(X,\Sigma,m)$ onto
 its linear subspace spanned by functions $h_1,\dots,h_k$.

\begin{theorem}
\label{th:1}
Suppose that the substochastic operator $P$ has the property {\rm (K)}.
Then there exist at most countable 
set~$J$, a family 
$\{A_j\}_{j\in J }$
of  pairwise disjoint
 measurable subsets of~$X$ 
and the family of substochastic operators  
$\{S_j\}_{j\in J }$, 
satisfying the conditions: 
\begin{enumerate} 
\item[\rm (i)] $P^*\mathbf 1_{A_j}\ge\mathbf 1_{A_j}$ for $j\in J $;
 \item[\rm (ii)] for each  $j\in J $ the operator $S_j$ 
is periodic on the set  $A_j$;
 \item[\rm (iii)]  for each  $j\in J $  and for any function
 $f\in L^1(X,\Sigma,m)$ we have
\begin{equation} 
\label{th:ap}
 \lim_{n\to \infty}\|\mathbf 1_{A_j} P^n\!f-S^n_jf\|=0;
\end{equation}
\item[\rm (iv)]  if $Y=X\setminus \bigcup_{j\in J } A_j$, 
then for any function $f\in L^1(X,\Sigma,m)$
and any compact set
  $F\subset X$ we have
\begin{equation} 
\label{th:sw}
\lim_{n\to \infty} \int\limits_{F\cap Y} P^n\!f(x)\, m(dx)=0. 
 \end{equation}
 \end{enumerate}
\end{theorem}

\begin{remark}
\label{r:as-okres}
From (i) it follows that 
the operator $P$ can be restricted to the space $L^1(A_j)$ 
and $P$ is a stochastic operator on this space.
Also the operator $S_j$ restricted to the space $L^1(A_j)$ is a stochastic operator and this restriction is given by the formula
\begin{equation} 
\label{wz:naSj}
S_jf=\sum_{i=1}^k \alpha_i(f) h_{i+1},\quad \text{$\alpha_i(f)=\int_{B_i} f(x)\,m(dx)$,   $\,B_i=\supp h_i$,}
\end{equation}
where the functions $h_1,\dots,h_k$ are the densities 
appearing in the definition of periodic operator.
The operator $P$ satisfies \eqref{per-op2}, i.e. $Ph_i=h_{i+1}$.
\end{remark}

The following criterion of sweeping follows from Theorems~\ref{th:2} and~\ref{th:1}.
\begin{corollary}
\label{col-sw}
Assume that a substochastic (continuous or discrete-time) semigroup $\{P(t)\}_{t\ge 0}$ satisfies
condition {\rm (K)} and has no invariant densities.
Then $\{P(t)\}_{t\ge 0}$ is sweeping from compact sets.
\end{corollary}
Our next aim is to find conditions guaranteeing that a stochastic semigroup $\{P(t)\}_{t\ge 0}$ 
satisfies the {\it Foguel alternative},\index{Foguel alternative} i.e. it
is asymptotically stable or sweeping from all compact sets  \cite{LiM}. The following result is a step in this direction. 
\begin{corollary}
\label{Foguel-al}
If a continuous-time substochastic semigroup $\{P(t)\}_{t\ge 0}$  is irreducible  and satisfies
condition {\rm (K)}, then the Foguel alternative holds.
\end{corollary}
In the case of substochastic operator we have a slightly weaker version of this 
corollary. 
A substochastic operator $P$ is called 
\textit{asymptotically periodic}\index{operator asymptotically periodic} if there exists a periodic operator $S$ such that 
\[
\lim_{n\to\infty} \|P^nf-S^nf\|=0\quad \text{for $f\in D$}.
 \]
 \begin{corollary}
 \label{c:op-a-F}
 If a substochastic operator satisfies {\rm (K)} and is irreducible, then it is either asymptotically periodic or sweeping from compact sets.
 \end{corollary}

Now we apply Corollary~\ref{Foguel-al} to  Markov chains. Let $Q$ be the intensity matrix of a continuous-time Markov chain~$(\xi_t)_{t\ge 0}$
 on $\mathbb N$. 
Since the space $\mathbb N$ is discrete the 
stochastic semigroup $\{P(t)\}_{t\ge 0}$ 
related to the process $(\xi_t)_{t\ge 0}$ satisfies condition (K).
\begin{theorem}
\label{af-tw1-c}
 Let  $\{P(t)\}_{t\ge 0}$ be the  stochastic semigroup on $l^1$ 
 related to a continuous-time Markov chain $(\xi_t)_{t\ge 0}$ 
with intensity matrix $Q$.
Let us assume that the entries of the matrix $Q$ satisfy the
following condition
\begin{enumerate}
\item [\rm (I)]  for all $i$ and $j$  there exists a sequence
$(i_1,i_2,\dots,i_m)$ such that  
$i_1=i$, $i_m=j$ and 
$q_{i_{k}i_{k+1}}>0\,$ for $k=1,\dots,m-1$.
\end{enumerate} 
Then the semigroup  $\{P(t)\}_{t\ge 0}$ satisfies the Foguel
alternative:

\begin{enumerate}
\item[\rm (a)] if the semigroup $\{P(t)\}_{t\ge 0}$ has an invariant
density, then it is asymptotically stable,

\item[\rm (b)] if the semigroup $\{P(t)\}_{t\ge 0}$ has no invariant
density, then for every  $x\in l^1$ and $i\in \mathbb N$ we have
\begin{equation}
\label{prz12af}
\lim_{t\to\infty} (P(t)x)_i=0.
\end{equation}
\end{enumerate}
\end{theorem}

If a stochastic semigroup (continuous or discrete-time) overlaps supports,
then the statements of Theorems~\ref{th:2} and~\ref{th:1} can be significantly strengthened (see \cite{PR-cell-cycle} Proposition~1).

\begin{corollary}
\label{prop:FA-SO}
Assume that 
a stochastic semigroup
$\{P(t)\}_{t\ge 0}$  satisfies {\rm (K)} and overlaps supports. Then
this semigroup  is sweeping or   has an invariant density $f^*$ with a support $A$  and there exists
a positive linear functional $\alpha$  defined on $L^1(X,\Sigma,m)$
 such that
\begin{itemize}
 \item[\rm(i)]\quad  for every  $f\in L^1(X,\Sigma,m)$ we have
\begin{equation}
\label{th2:ap2}
 \lim_{t\to \infty}\|\mathbf 1_{A} P(t)f-\alpha(f)f^*\|=0,
\end{equation}
\item[\rm(ii)]\quad if $\,Y=X\setminus A$, then for every  $f\in L^1(X,\Sigma,m)$ and for every compact set $F$ we have
\begin{equation}
\label{th2:sw2}
\lim_{t\to \infty} \int\limits_{F\cap Y} P(t)f(x)\, m(dx)=0.
 \end{equation}
 \end{itemize}
In particular, if $\{P(t)\}_{t\ge 0}$  has an invariant density $f^*$ with the support $A$ and  $X\setminus A$ is a subset of a compact set,
then this semigroup is asymptotically stable.
\end{corollary}

In order to formulate a further corollary, we need to introduce two additional
conditions, which are weak versions of well-known concepts: irreducibility
and tightness of measures.
\begin{enumerate}
\item[(WI)] there exists a point $x_0\in X$ such that for each $\varepsilon >0$ and for each density $f$ we have 
\begin{equation} 
\label{k:WI}
\int\limits_{B(x_0,\varepsilon)} P(t)f(x)\,m(dx)>0\quad\textrm{for some $t>0$}.
\end{equation}
\end{enumerate}
\begin{enumerate}
\item[(WT)] There exists $\kappa>0$ such that
\begin{equation}
\label{k:T}
\sup\limits_{F\in \mathcal F}\limsup_{t\to\infty} \int_F P(t)f(x)\,m(dx)\ge \kappa
\end{equation}
for $f\in D_0$, where $D_0$ is a dense subset of $D$ and $\mathcal F$ is the family of all compact subsets of $X$.
\end{enumerate}

\begin{corollary}
\label{cor:st}
Let $\{P(t)\}_{t\ge0}$ be a continuous-time stochastic semigroup.
Assume that $\{P(t)\}_{t\ge0}$ satisfies conditions {\rm (K)},  {\rm (WI)}, and {\rm (WT)}.
Then the semigroup $\{P(t)\}_{t\ge0}$ is  asymptotically stable.
\end{corollary}

\begin{corollary}
\label{compact-as}
Suppose that $X$ is a compact space. Then a continuous-time stochastic semigroup which satisfies {\rm (K)} and {\rm (WI)} is 
asymptotically stable.
\end{corollary}

The proof of Corollary~\ref{cor:st} can be found in \cite{PR-SD}. The second corollary follows directly from the first one.
New results concerning positive operators on Banach lattices similar in spirit to Theorems~\ref{asym-th2} and~\ref{th:2} 
may be found in \cite{Gerlach-Gluck}.

An advantage of the formulation of some results here 
in the form of the Foguel alternative  
 is that in order to show asymptotic stability we do not need 
to prove  the existence of an invariant density. It is enough to check that  
the semigroup is not sweeping with respect to compact sets then, 
automatically, the semigroup   $\{P(t)\}_{t\ge 0}$ is asymptotically stable.
We can eliminate the sweeping by means of some method similar to 
that of  Lyapunov function called Hasminski\u{\i} function.\index{Hasminski\u{\i} function}

\begin{theorem}
\label{tw:funkcja-Hasminskiego}
Let $A$ be the generator of a stochastic semigroup
$\{P(t)\}_{t\ge 0}$, 
$V\colon X\to [0,\infty)$ be a measurable function,
$Z$ be a measurable set, and $\varepsilon>0$.  
Assume that there exists a measurable function, formally denoted by $LV$, satisfying conditions
\begin{equation}
\label{c:Hf-i1}
\int_X LV(x) g(x)\,m(dx)=\int_X  V(x)Ag(x)\,m(dx)\quad\text{for $g\in \mathfrak D(A)\cap D_V$},
\end{equation}
\begin{equation}
\label{c:Hf-i2}
 LV(x)\le M-\varepsilon \quad\text{for $\,x\in Z$}, \quad\quad LV(x)\le -\varepsilon \quad\text{for $\,x\in X\setminus Z$},
\end{equation}
where $D_V=\{f\in D\colon \,\int_Xf(x)V(x)\,m(dx)<\infty\}$.
Then the semigroup $\{P(t)\}_{t\ge 0}$ is not sweeping with respect to the set  $Z$. If $Z$ is compact, then  the semigroup satisfies condition $(\rm{WT})$.
\end{theorem}

The main difficulty in applying this method lies in correctly defining the action of the operator $L$ on the function $V$, since $V$ does not usually belong to the domain of the operator $L$.
The Hasminski\u{\i} function method has been used to study Markov chains~\cite{PRT2012}, multistate diffusion processes~\cite{PR-jmaa1}, diffusion with jumps~\cite{PR-bpan}, dynamical systems with jumps~\cite{Pichor1998}. 

\subsection{Constrictive operators}
\label{ss:constrictive}
The next theorem we discuss here shows that convergence of trajectories of a stochastic semigroup to a compact set leads to interesting asymptotic properties. 
A stochastic semigroup $\{P(t)\}_{t\ge 0}$
is called  strongly (weakly) \textit{constrictive}\index{constrictive semigroup} \cite{LLY}
if there exists a strongly (weakly) compact subset $F$ of
$L^1$ such that
\[
\lim_{t\to\infty}\inf\{\|P(t)f-g\|\colon g\in F\}=0
\quad\text{for $f\in D$}.
\]

Even though the definition of the constrictive operator and the property (K) are significantly differ, similar conclusions follow. One of them is the following theorem. 

\begin{theorem}
\label{tw:Komornika} 
Suppose that a stochastic operator $P$ defined on the space 
$L^1(X,\Sigma,m)$  
is constrictive.
Then there exist a finite set~$J$, a family
$\{A_j\}_{j\in J}$
 of pairwise disjoint measurable subsets of~$X$ 
and the family
$\{S_j\}_{j\in J}$
of substochastic operators
satisfying the conditions: 
\begin{enumerate} 
\item[(i)] $P^*\mathbf 1_{A_j}\ge \mathbf 1_{A_j}$ for $j\in J $;
 \item[(ii)] for each  $j\in J $ the operator $S_j$ is periodic 
on the set $A_j$;
 \item[(iii)] for each  $j\in J $ and for any function $f\in L^1(X,\Sigma,m)$ we have
\begin{equation} 
\label{th:ap-Kom}
 \lim_{n\to \infty}\|\mathbf 1_{A_j} P^n\!f-S^n_jf\|=0;
\end{equation}
\item[(iv)]  if $Y=X\setminus \bigcup\limits _{j\in J} A_j$, then
for any function $f\in L^1(X,\Sigma,m)$
we have
\begin{equation} 
\label{th:sw_Kom}
\lim_{n\to \infty} \int\limits_{Y} P^n\!f(x)\, m(dx)=0. 
 \end{equation}
 \end{enumerate}
\end{theorem}
Theorem~\ref{tw:Komornika}  was first proved in~\cite{LLY} under  assumption of strong constrictivity; two years later Komorn\'{\i}k\index{Komorn\'{\i}k} \cite{Komornik} showed that weak constrictivity suffices.  
The original formulation of Theorem~\ref{tw:Komornika}  was different from that cited here. We have given its statement formulated in the spirit of Theorems~\ref{th:2} and~\ref{th:1}.

If there exists a function $h\in L^1$ such that 
\begin{equation}
\label{ufu}
\lim_{t\to\infty}\|(P(t)f-h)^+\|=0\quad\textrm{for every $f\in D$},
\end{equation} 
then Theorem~\ref{tw:Komornika} holds. It is a simple consequence of the fact that the set $F= \{f\in D\colon f\le h\}$ is weakly compact.
Theorem~\ref{tw:Komornika} can be strengthened considerably in the case when $\{P(t)\}_{t\ge 0}$ is a continuous-time stochastic semigroup. 
Then we replace condition~\eqref{th:sw_Kom} by condition~\eqref{pomoc:ogolny}.

The property of constrictivness was generalised and studied in many papers (see \cite{Emelyanov} for a comprehensive overview of these results).
Examples of applications of Theorem~\ref{tw:Komornika}
and the lower function theorem
are given in \cite{LiM}.

\section{Applications to biological models}
\label{s:biol-appl}
\subsection{Model of DNA evolution}
\label{ss:DNA-ewolucja}
The degree of  relatedness between humans, species, genes and proteins, can be studied by comparing DNA chain sequences.
DNA is inherited from parents, but mutations occur as a result of errors in DNA replication and repair.
Mutations occur throughout the DNA chain at a constant level on average, and although the potential
mutations may be numerous, as a result of natural selection, most of them are not conserved and
mainly selection-neutral mutations are observed.

The simplest way to measure evolutionary distance is to check in how many places the aligned DNA sequences differ.
More precisely, we first compare and align DNA sequences so that they match, and then check the differences between them. 
The initial issue here is modelling nucleotide substitutions at a fixed point in the DNA chain.
Since we have four different nucleotides corresponding to nucleobases A, G, C, T, so models describing nucleotide evolution, 
called \textit{models of substitution} of nucleotides in the
DNA sequences\index{substitution model} 
are Markov chains on a four-state space
$S=\{{\rm A,G,C,T}\}$. The intensity matrix in this chain is of the form
\[
Q=
\left[
\begin{matrix}
q_{\rm AA}&q_{\rm AG}&q_{\rm AC}&q_{\rm AT}\\
q_{\rm GA}&q_{\rm GG}&q_{\rm GC}&q_{\rm GT}\\
q_{\rm CA}&q_{\rm CG}&q_{\rm CC}&q_{\rm CT}\\
q_{\rm TA}&q_{\rm TG}&q_{\rm TC}&q_{\rm TT}
\end{matrix}
\right].
\]
 
The simplest substitution model is 
\textit{Jukes--Cantor model}\index{Jukes--Cantor model} \cite{JC}.
 In this model, we assume that
probability of replacing any nucleotide, by another nucleotide is the same. The intensity matrix $Q$ in this model is of the form
\[
Q=
\left[
\begin{matrix}
-3\lambda &\lambda&\lambda&\lambda\\
\lambda&-3\lambda&\lambda&\lambda\\
\lambda&\lambda&-3\lambda&\lambda\\
\lambda&\lambda&\lambda&-3\lambda
\end{matrix}
\right].
\]
Let $p_{ij}(t)={\rm P}(\xi_t=j|\xi_0=i)$. 
It is easy to check that
\[
p_{ij}(t)=\frac14-\frac14e^{-4\lambda t} \text{ \,\,for $i\ne j$,}\quad\  
p_{ii}(t)=\frac14+\frac34e^{-4\lambda t} \text{ \,\,for $i\in S$.}
\]

These formulas allow us to estimate the evolutionary distance
 of two different individuals or species based on differences in aligned DNA sequences.
We will assume that both sequences are derived from a common ancestor and that, as a result of mutations, there have been independent changes in both sequences in successive generations.
As a measure of evolutionary distance, we will take the time elapsed since the common ancestor, or the number of generations.

We assume that a mutation of a single nucleotide occurred on average every $t_0$ years. To begin, we compare the sequence of an ancestor that lived $t$ years ago with that of a contemporary living individual.   
Since the intensity of nucleotide change is $3\lambda$,  so $3\lambda t_0=1$.
The probability of changing a nucleotide at a fixed site is
\[
p= 3\bigl(\tfrac 14-\tfrac14 e^{-4\lambda t}\bigr).
\]
Since $3\lambda t_0=1$, we have 
\begin{equation}
\label{dystans-JC}
\frac{t}{t_0}= -\tfrac{3}{4}\ln\big(1-\tfrac 43 p\big).
\end{equation}
If the sequences are of length $n$ and we have $n_d$ differences, we assume that 
$\hat p=n_d/n$ is the estimated value of $p$.

Let us now consider two sequences of contemporary individuals. Suppose we want to determine how many years ago their common ancestor lived. 
The Markov chain in this model is reversible.
Thanks to this symmetry, we can treat one of the individuals as the ,,ancestor'' of their common ancestor. Thus, the time to the common ancestor will be twice as small as the determined by the formula~(\ref{dystans-JC}).

\subsection{Growth of erythrocyte population}
\label{ss:z-erythrocyte}
Consider the following 
\textit{model of erythrocyte population growth}\index{model of erythrocyte population growth} (red blood cells).
Erythrocytes are formed in the bone marrow by the division and differentiation of erythroid stem cells. Mature erythrocytes enter the bloodstream, where they serve
mainly to transport oxygen to the cells and are destroyed after about four months. Assume that $b>0$ is the intensity of erythrocyte production in the bone marrow,
and $d>0$ is the destruction rate of a single erythrocyte. We check that the number of erythrocytes in the bloodstream is described 
by a Markov chain (a certain birth-death process)
on the space $\mathbb N$. The intensity matrix $Q$ has entries:
$q_{i\,i+1}=b$, $q_{i\,i-1}=di$, and $q_{ii}=-(b+di)$  for $i\in \mathbb N$ and $q_{ij}=0$ otherwise. 

First, using Theorem~\ref{Kol} we prove that the matrix $Q$ generates a stochastic semigroup $\{P(t)\}_{t\ge 0}$. 
Let  $\lambda>0$ be a positive constant and assume a contrary that there exists a nonzero and nonnegative $x\in l^{\infty}$ which is a  solution  
of the equation $Qx=\lambda x$. Then
\begin{equation}
\label{e:er1}
dix_{i-1}-(b+di)x_i+bx_{i+1}=\lambda x_i\quad \text{for $\,i\in \mathbb N$},
\end{equation}
where we formally assume that $x_{-1}=0$.
Then 
\begin{equation}
\label{e:er2}
x_{i+1}- x_i=\frac{\lambda}{b} x_i+\frac db(x_i-x_{i-1})i.
\end{equation}
Hence, the sequence $(x_i)$  is increasing and $x_{i+1}\ge \big(1+\frac{\lambda}{b}\big) x_i$.
Consequently, $x\notin l^{\infty}$,  and the matrix $Q$ generates a stochastic semigroup.

We show that the semigroup $\{P(t)\}$ has an invariant density.
Since the generator $A$ of this semigroup is given by the formula $Ax=xQ$, we are looking for a positive solution of the equation $xQ=0$.
Then
\begin{equation}
\label{e:er3}
bx_{i-1}-(b+di)x_i+d(i+1)x_{i+1}=0\quad \text{for $\,i\in \mathbb N$},
\end{equation}
where $x_{-1}=0$. The  sequence $x_i=b^i/(d^ii!)$ 
satisfies \eqref{e:er3}.
Since $\sum_{i=0}^{\infty}x_i=e^{b/d}$,   
in order to obtain a probabilistic distribution we divide $x_i$ by $e^{b/d}$. This will yield the invariant density $(\pi_i)$,
which is a Poisson distribution with the parameter $\lambda=b/d$:
\[
\pi_i=\frac{\lambda^ie^{-\lambda}}{i!} \quad \text{for $\,i\in \mathbb N$.}
\]
Since the matrix $Q$ satisfies condition (I) of Theorem~\ref{af-tw1-c},
the semigroup $\{P(t)\}_{t\ge 0}$ is asymptotically stable.

\subsection{Population growth with stochastic noise}
\label{ss:z-dyfuzji}
Consider the following equation 
\begin{equation}
\label{rs-pop1}
dx_t=b(x_t)\, dt+\sigma(x_t)\,dw_t. 
\end{equation}
Eq.~\eqref{rs-pop1} is a stochastic version of the deterministic models describing the growth  of a single population $x'(t)=b(x)$, where 
$x(t)$ is the size of the population at time $t$ and 
the function $b$ is the growth coefficient. We assume that  
$b\colon [0,\infty)\to\mathbb R$ 
is a $C^1$-function
on the whole interval $[0,\infty)$ satisfying additionally the conditions $b(0)=0$ and $b(x)\le Kx$ for $x\ge 0$ and some $K>0$.

We can consider two different stochastic perturbations. 
The first, called \textit{environmental noise},\index{environmental noise}
is used to describe small populations in which random fluctuations act in the same way on all individuals in the population, so $\sigma(x)=\sigma x$, $\sigma>0$. 
The second is called \textit{demographic noise}.\index{demographic noise} 
In this case the stochastic perturbation 
acts independently on each part of population
and it has zero mean. According to the central limit theorem 
such noise is proportional to 
$\sqrt {x}$, so $\sigma(x)=\sigma \sqrt {x}$.   
Models with demographic noise are more difficult to study because the trajectories of the process $(x_t)$ 
 are absorbed with positive probability at $0$.  

We only consider here models with environmental noise, thus Eq.~\eqref{rs-pop1} has the form 
\begin{equation}
\label{rs-pop2}
dx_t=b(x_t)\, dt+\sigma x_t \,dw_t. 
\end{equation}
The stochastic semigroup $\{P(t)\}_{t\ge 0}$ related to the process 
$(x_t)_{t\ge 0}$ is integral and the kernel $q(t,x,y)$ is positive and continuous for $t,x,y>0$. 
Thus, the semigroup $\{P(t)\}_{t\ge 0}$ is irreducible and satisfies 
condition~(K).
Then 
 \begin{equation}  
\label{FP-one}
\frac{\partial u}{\partial t}=
\frac12\frac{\partial^2}{\partial x^2}(\sigma^2x^2u)-
\frac{\partial}{\partial x}(b(x)u)
\end{equation}
is the Fokker--Planck equation corresponding to this semigroup.
An invariant density $f^*$, if it exists, 
is a solution of the equation
\begin{equation}  
\label{FP-one-st}
\frac12 \big (\sigma^2x^2f^*(x)\big )''-\big (b(x)f^*(x)\big )'=0.
\end{equation}
We easily determine from~\eqref{FP-one-st}  the invariant density 
\begin{equation}  
\label{FP-stat}
f^*(x)=\frac{C}{x^2}\exp\bigg(\int_{1}^x\frac{2b(s)}{\sigma^2s^2}\,ds\bigg),
\end{equation}
as long as $C$ can be chosen such that $f^*$ is a density.
Thus, an invariant density exists, and consequently the semigroup 
$\{P(t)\}_{t\ge 0}$ is asymptotically stable 
(see Corollary~\ref {asym-irreducible}),  if and only if 
\begin{equation}  
\label{FP-stat-war}
\int_0^{\infty}
\frac{1}{x^2}\exp\bigg(\int_{1}^x\frac{2b(s)}{\sigma^2s^2}\,ds\bigg)
\,dx
<\infty.
\end{equation}

If \eqref{FP-stat-war} does not hold then according to Corollary~\ref{col-sw}
the semigroup  $\{P(t)\}_{t\ge 0}$ is sweeping from intervals 
$[\alpha,\beta]$, $0<\alpha<\beta<\infty$, i.e. 
\begin{equation}  
\label{FP-stat-sweep}
\lim_{t\to\infty}
\int_{\alpha}^{\beta} P(t)f(x)\,\,dx=0\quad \text{for  $f\in D$.}
\end{equation}

Using Theorem~\ref{KT}, we can show more precise sweeping conditions and give their biological interpretation. If \eqref{FP-stat-war} does not hold,
then $f^*$ given by \eqref{FP-stat} with $C>0$ is 
a non-integrable positive  function such that $P(t)f^*=f^*$.
Assume that there exists the limit $b'(\infty)=\lim_{x\to\infty}b(x)/x$.
It is easy to check that 
\begin{align}
\label{wym-dp1}
&\int_0^{\alpha}f^*(x)\,dx<\infty\textrm{ if $\sigma^2<2b'(0)$},
 &&\int_0^{\alpha}f^*(x)\,dx=\infty\textrm{ if $\sigma^2>2b'(0)$},
\\
\label{wym-dp2}
&\int_{\beta}^{\infty}f^*(x)\,dx<\infty\textrm{ if $\sigma^2>2b'(\infty)$},
&&\int_{\beta}^{\infty}f^*(x)\,dx=\infty\textrm{ if $\sigma^2<2b'(\infty)$}
\end{align}
for $0<\alpha<\infty$ and $0<\beta <\infty$. 
From Theorem~\ref{KT} it follows that:
\begin{enumerate}
\item[(i)] if $\sigma^2<\min \{2b'(0),2b'(\infty)\}$, then 
$\lim\limits_{t\to\infty} {\rm P}(x_t\ge \beta)=1$,
\item[(ii)] if $\sigma^2>\max \{2b'(0),2b'(\infty)\}$, then 
$\lim\limits_{t\to\infty} {\rm P}(x_t\le \alpha)=1$,
\item[(iii)] if $2b'(0)<\sigma^2<2b'(\infty)$, then 
$\limsup\limits_{t\to\infty} {\rm P}(x_t\le \alpha)>0$
and
$\limsup\limits_{t\to\infty} {\rm P}(x_t\ge \beta)>0$.
\end{enumerate}
It can be check that the limits superior in (iii) can be replaced by limits, but  
their values depend on the distribution of $x_0$.
Let $b_0:= b'(0)-\sigma^2/2$ and $b_{\infty}:= b'(\infty)-\sigma^2/2$.
Then $b_0$ and $b_{\infty}$ are called 
the \textit{actual growth rate}\index{actual growth rate} at $0$ and $\infty$.
Now conditions (i)--(iii) say:
if both actual growth rates are positive, then the population grows indefinitely; if both actual growth rates are negative, 
then the population is becoming extinct; if $b_0<0$ and $b_{\infty}>0$, then the population can with positive probability both grow indefinitely and
become extinct.

\subsection{Gene expression}
\label{ss:gene-exp}
Gene expression is a process of production and decay of various biomolecules dependent on activity of a number of genes. 
Changes in  concentration of these biomolecules can be basically described by 
a  system of differential equations, but this system depends on random parameters related to the genes' activity. 

Assume that genes can form $k+1$ configurations depending on their state of activity
and for configuration $i\in I=\{0,\dots,k\}$
concentration of $d$ different
 biomolecules changes according to a system of differential equations of the form $x'=b^i(x)$,  $x\in \mathbb R_+^d$. 
We assume that for each $x_0\in \mathbb R_+^d$ there exists a unique solution 
$x\colon \mathbb R_+ \to \mathbb R_+^d$ of the last equation with the initial condition $x(0)=x_0$ and we denote this solution by $\pi_t^ix_0$.
Assume that for $i,j\in I$ and  $i\ne j$ and for concentration 
$x$, the rate of switching configuration from $i$ to $j$ 
is a nonnegative continuous and bounded function  
$x\mapsto q_{ij}(x)$.  

We start with a simple model in which
protein  molecules are directly "produced" by a single gene.
A gene can be in an active $i(t)=1$ or inactive state $i(t)=0$.
 When gene is active then proteins are produced with a constant speed $P$.
In both states $\mu$ is the rate of protein degradation.
Let $x(t)$ be the concentration of protein. 
Then $x(t)$ satisfies the equation
\begin{equation}
\label{1:ga1}
x'(t)=P i(t) -\mu x(t).
\end{equation}
Gene expression is described by
dynamical systems with random switching
 $\xi_t=(x(t),i(t))$, $t\ge 0$,
with values in  the space $\mathbb R_+\times \{0,1\}$
and
with random jumps from 
$(x,0)$ to $(x,1)$ and from
$(x,1)$ to $(x,0)$ with rates $q_{01}(x)$ and
$q_{10}(x)$, respectively.
According to Sec.~\ref{ss:diff-with-switching} this process generates 
a stochastic semigroup on the space $L^1(X,\mathcal B(X),m)$, where $X=\mathbb R_+\times \{0,1\}$
and $m$ is the product of the Lebesgue measure on $\mathbb R_+$ and the counting measure on the set $\{0,1\}$.

\begin{figure}
\begin{center}
\begin{picture}(280,100)(0,0)
\includegraphics[scale=1.0,viewport=142 600 456 364]{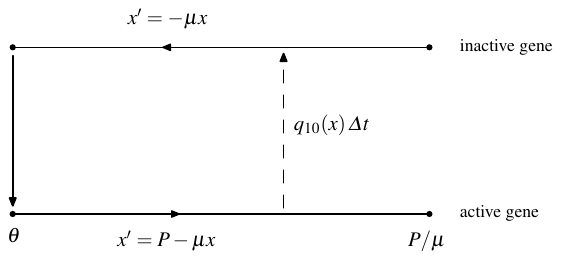}
\end{picture}
\end{center}
\caption{Gene expression scheme with  jumps on the boundary.}
\label{r:genI-ze-skok}
\end{figure}

We now present a version  of this  model with non-random moments of jumps. 
We assume that the concentration  
of protein must not fall below a fixed level $\theta>0$.
We assume that the gene is inactivated with rate $q_{10}(x)$ but
if the protein concentration reaches the level $\theta$, the  
protein production is automatically activated (Fig.~\ref{r:genI-ze-skok}).

Now gene expression is described by
a piecewise deterministic Markov process
 $\xi_t=(x(t),i(t))$, $t\ge 0$,
with values in  the space 
$\big[\theta,\frac P{\mu}\big]\times \{0,1\}$
with random jumps form $(x,1)$ to $(x,0)$ and
with non-random jumps from $(\theta,0)$ to $(\theta,1)$.

If $x_0$ the initial number of protein in the active state, then the length of the active state has the distribution 
\begin{equation}
\label{d:prost-ge}
F_1(t)=1-\exp\bigg\{- \int_0^t q_{10}(1 +(x_0-1)e^{- s}) \,ds\bigg\}.
\end{equation}
The duration  of stay in the inactive state $T$ satisfies the equation
$x_0e^{-\mu t}=\theta$,  thus 
$T=\frac1{\mu}\ln(x_0/\theta)$. 
The length distribution of the inactive state
is therefore the function
\begin{equation}
\label{d:o}
F_0(t)=\begin{cases}
0, &\text{if $t< T$},\\
1, &\text{if $t\ge  T$}.
\end{cases}
\end{equation}

Now we consider more realistic  model introduced in \cite{lipniacki2} and studied in \cite{BLPR}.
A model consists of two types of  molecules:
mRNA, and protein, and we assume that variables $x_1,x_2$ denote
their concentrations.
If the gene is active then
mRNA transcript molecules are synthesized
with a constant speed $R$.
The protein translation proceeds with the rate
$Px_1(t)$, where $P$ is a constant.
The mRNA and protein
degradation rates are $\mu_R$ and $\mu_P$,
respectively.
Thus the process is described by the following system
\begin{equation}
\label{1:ga1b}
\left\{
\begin{aligned}
x_1'(t)&=Ri(t) -\mu_R x_1(t),\\
x_2'(t)&=Px_1(t) -\mu_P x_2(t),
\end{aligned}
\right.
\end{equation}
where $i(t)=1$ if the gene is active and
$i(t)=0$ in the opposite case.

We assume that the rates $q_{01}$ and $q_{10} $  are continuous and bounded functions of  $\mathbf x=(x_1,x_2)$.
Then the process $\xi_t=(x_1(t),x_2(t),i(t))$, $t\ge 0$,
is given by dynamical systems with random switching.
The trajectories of the process enter the invariant set
$X=\Big[0,\frac{R}{\mu_R}\Big]\times \Big[0,\frac{PR}{\mu_P\mu_R}\Big] \times \{0,1\}$, so we can consider the process  $(\xi_t)_{t\ge 0}$ with values in this set. 
It generates a stochastic semigroup 
$\{P(t)\}_{t\ge 0}$
on the space 
$L^1(X,\mathcal B(X),m)$.

In this model we have two dynamical systems,  which correspond, respectively, to $i=0$ and $i=1$.
Each system  has a unique stationary point, either $\mathbf x_0=(0,0)$, or $\mathbf x_1=\Big(\frac{R}{\mu_R},\frac{PR}{\mu_P\mu_R}\Big)$, which is asymptotically stable.
We assume that
$q_{01}(\mathbf x_0)>0$ and $q_{10}(\mathbf x_1)>0$.
This assumption implies that the trajectories of the process  visit any neighbourhood of the point 
$(\mathbf x_0,0)$ and the  point
$(\mathbf x_1,1)$ for infinitely many times. 
It means that condition (WI) holds. Moreover, if we assume additionally that
$q_{01}(\mathbf x_1)>0$ or  $q_{10}(\mathbf x_0)>0$ then we can check that condition (K) also holds (see \cite{PR-prz-apss-apdMp} for the proof based on 
H{\"o}rmander's condition).
The set $X$  is compact. 
According to Corollary~\ref{compact-as}, 
the semigroup
$\{P(t)\}_{t\ge0}$ is asymptotically stable.

In some cases, the process of protein production has 
a pre-transcriptional step. 
Activation of the gene is followed by the production of RNA molecules called pre-mRNA,
which cannot directly produce a protein. 
Then, while still in the cell nucleus, non-coding sequences (introns)
are truncated to form mRNA molecules, which are then transported to the cytoplasm, where 
they serve as matrices in the translation  process of 
protein production. 
Such a three-step model was introduced and investigated in~\cite{RT}.

Let us denote by $x_1(t),x_2(t),x_3(t)$ the  concentration of pre-mRNA, mRNA and protein at time $t$. 
The change in these quantities over time can be represented by a system of three equations
\begin{equation}
\label{1:ga1b-2}
\left\{
\begin{aligned}
x_1'(t)&=Ai(t) -(R+\mu_{pR}) x_1(t),\\
x_2'(t)&=R x_1(t) -\mu_R x_2(t),\\
x_3'(t)&=P x_2(t) -\mu_P x_3(t),
\end{aligned}
\right.
\end{equation}
where $A$ is the speed of pre-mRNA molecules synthesis  when the gene is active,
$R$ is the conversion rate of pre-mRNA molecules into active mRNA,
$\mu_{pR}$ is the degradation rate of pre-mRNA molecules, and the other constants are the same as in the previous model.

The model \eqref{1:ga1b} can also be modified to a model with non-random jumps.
As in the previous model the gene is inactivated with rate $q_{10}(x_1,x_2)$
but it is activated when the protein concentration $x_2(t)$ reaches $\theta$,
where $0<\theta\le PR/\mu_R\mu_P$.
Note that in this case it may happen that 
$x_1(\bar t)<\mu_p\theta $ at the time $\bar t$ of activation, 
and then $x_2'(\bar t)<0$ and the protein concentration $x_2(t)$ 
will decrease in a certain interval
$(\bar t,\bar t+\varepsilon)$,
in particular will fall below the level of $\theta$.
It is clear that the gene should remain active until the concentration
of the protein $x_2(t)$ does not exceed at least the level of $\theta$,
 we should therefore assume that the inactivation intensity
$q_{10}(x_1,x_2)$ is zero when $x_2\le \theta$.
The  process 
$\xi_t=(x_1(t),x_2(t),i(t))$, $t\ge 0$
is defined on the space 
$X=\big[0,\frac R{\mu_R}\big]\times \big[0,\frac {PR}{\mu_P\mu_R}\big]\times \{0,1\}$ and has non-random jumps on the set 
$S=\big[0,\frac R{\mu_R}\big]\times \{\theta\}
\times \{0,1\}$.
In this case, additional jumps occur inside the phase space.

The papers \cite{HWS,KKLL} present more advanced models 
describing the process of 
 \textit{subtiline} production \index{subtiline production}
 by the bacterium   \textit{Bachillus subtilis}.
  Although the process plays an important role in biochemical technology 
 we will not discuss the model in detail, but only highlight its
interesting mathematical aspects.
 
 To survive, the bacterium produces a type of antibiotic,  subtiline,
which eliminates other microorganisms competing with it.
The model considered in \cite{HWS} consists of five differential equations 
describing the size $x_1$ of the \textit{B. subtilis} population, 
the amount of food $x_2$ and the concentration $x_3,x_4,x_5$ of the three proteins involved in the production.
Two of these equations are without switching, two with switching 
dependent on the state of two genes and one with switching dependent on the amount of food:
\[
x_3'(t)=
\begin{cases}
-\mu x_3(t), &\text{if $x_2(t)\ge \theta$},\\
P_3-\mu x_3(t), &\text{if $x_2(t)<\theta$},
\end{cases}
\]
where $\theta$ is a positive constant.
The corresponding process is defined on the space $\mathbb R^5_+$
and has additional jumps on the subspace 
$S=\{x\in \mathbb R^5_+\colon x_2=\theta\}$.

\subsection{Movement with  velocity jumps}
\label{ss:transport}
The pure jump-type processes considered in Sec.~\ref{ss:kanngaroo-movement}
are a far-reaching simplification of the description of the movement of biological objects. 
Much closer to reality models of the movement of cells or organisms are 
\textit{velocity jump processes}~\cite{ODA}.\index{velocity jump process}    
These processes are also used to describe vesicular transport in cells.
We will start with a general description and then give two examples of applications.

An individual is moving in the space $\mathbb R^d$ with a constant velocity and at jump times  $(t_n)$
it chooses a new velocity. 
The intensity of the jump $\psi$ depends on its location $x\in \mathbb R^d$ and velocity $v\in V$,
where $V=\mathbb R^d$ or $V$ is some subset of $\mathbb R^d$, e.g. $V$ is the unit sphere $S^{d-1}$ in $\mathbb R^d$.
Assume that $\psi\colon \mathbb R^d\times V\to [0,\infty)$ 
is a bounded measurable function. 
Let $x,v$ be the position and velocity  at jump time. Denote by $T$ the time until the next jump.
Then 
\[
{\mathrm P}(t<T\le t+\Delta t \mid T>t)=\psi(x+vt,v)\Delta t+o(\Delta t).
\]
We can write this formula using the cumulative distribution function $F$ of the random variable $T$:
\[
\frac{F(t+\Delta t)-F(t)}{1-F(t)}=\psi(x+vt,v)\Delta t+o(\Delta t).
\]
Dividing both sides by $\Delta t$ and going to the limit when
$\Delta t\to 0$
we obtain
\[
\frac{F'(t)}{1-F(t)}=\psi(x+vt,v).
\]
We write the last equation  as
\[
\frac{d}{dt} \ln(1-F(t))=-\psi(x+vt,v),
\]
and then we solve it and obtain
\[
F(t)=1-e^{-\int_0^t\psi(x+vs,v)\,ds}. 
\]

Let $x(t)$  be the position and $v(t)$ be the velocity of an individual at time $t$.
 We assume that for every $x\in \mathbb R^d$, $v\in V$,
there is 
a probability Borel measure $\mathcal P(x,v,B)$ on $V$ which describes
the change of the velocity after a jump, i.e.,
\[
{\rm P}(v(t_n)\in B\,\mid x(t_n^{-})=x,\, v(t_n^{-})=v)=\mathcal P(x,v,B)
\]
for every Borel subset $B$ of $V$, where
$x(t_n^{-})$ and $v(t_n^{-})$ are the left-hand side limits of
$x(t)$ and $v(t)$, respectively, at the point $t_n$.

Knowing the distribution of the time to the next jump and the measure $\mathcal P(x,v,B)$ 
describing the change in velocity after the jump, we can construct
a piecewise homogeneous deterministic Markov process $(\xi_t)_{t\ge 0}$
on the space $X=\mathbb R^d\times V$.
Such processes have been used in~\cite{BB,ODA,Stroock}  to describe the movement of bacteria, insects and mammals.

Let $m$ be a $\sigma$-finite measure on $(V,\mathcal B(V))$.
We usually take it to be the Lebesgue measure when the set $V$ is of positive Lebesgue measure, or
the Lebesgue-Riemann surface measure when $V$ is a manifold. 
Suppose that, for every point $x\in\mathbb R^d$  the transition probability function $\mathcal P(x,v,B)$ 
corresponds to some stochastic operator $P_x$ on $L^1(V,\mathcal B(V),m)$
and $\psi$ is a bounded function. Then the operator
\begin{equation}
\label{op-infi-skok-veloc-jump}
Af(x,v)= -v\cdot \nabla_x f(x,v)+ P_x(\psi f)(v) - \psi(x,v)f(x,v)
\end{equation}
is the generator of a stochastic semigroup $\{P(t)\}_{t\ge 0}$ on 
$L^1(X,\mathcal B(X),dx\times m)$. 
As in Sec.~\ref{ss:kanngaroo-movement}, the operator $A$ can be written in the form 
\begin{equation}
\label{op-infi-skok-veloc-jump2}
Af=-v\cdot \nabla_x f +\lambda \bar P_xf-\lambda f.  
\end{equation}

\begin{example}
\label{telegraph}
Consider a  symmetric movement on the real line $\mathbb R$.
In this case we assume that an individual is  moving with
 constant speed, say one, and at a jump time
it changes the direction of movement to the opposite one.
The process $(\xi_t)_{t\ge 0}$
 has values in the space $X=\mathbb R\times\{-1,1\}$
and  $\mathcal P(x,v,\{-v\})=1$ for $v=-1,1$.
We assume that 
jump times
are determined by a Poisson process $(N_t)_{t\ge 0}$ with intensity $\lambda$, it means that $\psi\equiv \lambda$.
We can give explicit formulae for the  process $(\xi_t)_{t\ge 0}$:
\begin{equation}
\label{r:Kac-bezp}
v_t=v_0(-1)^{N_t},\quad x_t=x_0+v_0\int_0^t(-1)^{N_s}\,ds.
\end{equation}
The stochastic semigroup related to this process is defined on the space $L^1(X,\mathcal B(X),m)$, where 
 $m$ is the product of the Lebesgue measure on $\mathbb R$ and 
the counting measure on the set $\{-1,1\}$.
If the random variable  $\xi_0=(x_0,v_0)$
has a density $u_0(x,i)$, then the random variable
$\xi_t=(x_t,v_t)$
has the density $u(t,x,i)$ satisfying the following system of equations:
\begin{equation}
\label{r:Kaca}
\begin{aligned}
\frac{\partial u(t,x,1)}{\partial t}&=-\frac{\partial u(t,x,1)}{\partial x}-\lambda u(t,x,1)+\lambda u(t,x,-1),\\
\frac{\partial u(t,x,-1)}{\partial t}&=\frac{\partial u(t,x,-1)}{\partial x} -\lambda u(t,x,-1)+\lambda u(t,x,1),\\
u(0,x,i)&=u_0(x,i).
\end{aligned}
\end{equation}
We recall that the function  
$u(t,x,i)$ is a density of $\xi_t$ if
\begin{equation*}
{\rm P}(\xi_t\in A\times \{1\}\cup B\times \{-1\})=\int_A u(t,x,1)\,dx+\int_B u(t,x,-1)\,dx,
\end{equation*}
where $A$ and  $B$ are any Borel subsets of $\mathbb R$.
The system (\ref{r:Kaca}) has been studied by Goldstein \cite{Goldstein} and Kac \cite{Kac-r-walk}, among others.
By suitable substitution it can be reduced to the telegraph equation.
\end{example}

There are many interesting examples of velocity jump processes
 with applications to aggregation and chemotaxis
phenomena (see e.g. \cite{HH}). 
Stochastic billiards introduced in Sec.~\ref{ss:proc-non-random-jumps}
is also a velocity jump process but it has jumps at the boundary,
so its description is different.
Also the cell cycle can be  modeled 
by a velocity jump process~\cite{Rotenberg,LeR}.

Other interesting velocity jump processes
with jumps on certain subsets are used in  
\cite{BN1,BN2} to describe vesicular transport in cells.
In these models, a molecule moves along the interval $[0,L]$ to some destination point $x_0$. 
A molecule can be in three states
$i\in I=\{-1,0,1\}$, we assume that in the state $i$ 
it moves with a velocity $i$. The transition between states is controlled by a homogeneous Markov chain 
$\vartheta(t)$ on the set $I$.
The particle starts from a point $x=0$ and moves to the right.
Assume that $0$ is an elastic screen and $L$ is an absorbing screen.
We assume that if a particle is in a certain neighborhood $U$ of the point $x_0$ and in the state $i=0$, it is 
absorbed by the target point with intensity $\kappa$.
If the molecule is captured by the target point or the point $L$, 
the process repeats itself.
Between jumps the molecule moves
according to a jump process $\xi_t=(x(t),i(t))$, $t\ge 0$, 
defined on the phase space $[0,L]\times I$,
where the pair $(x(t),i(t))$ satisfies the system of equations
\begin{equation}
\left\{
\begin{aligned}
x'(t)&=i(t),\\
 i'(t)&=0.
\end{aligned}
\right.
\end{equation}
If $x\in (0,L)$, then the process $\xi_t$ can jump from $(x,i)$ to
$(x,j)$ with intensity $q_{ij}$ for $i,j\in I$;
if $(x,j)\in U\times \{0\}$, the process can jump from $(x,0)$ to
$(0,1)$ with intensity~$\kappa$.
The process has additional jumps from $(0,-1)$ and from $(L,1)$ to $(0,1)$
with probability~$1$.

The literature on stochastic models of intracellular transport is quite extensive.
Interested readers are referred to the review paper~\cite{BN3}.

\subsection{Structured population models}
\label{ss:structured}
A population is usually heterogeneous and it is important to divide
it into homogeneous groups according to some significant parameters
such as age, size, maturity or trait in phenotypic models.
Models which describe the time evolution of the distribution of the population
according to such set of fixed parameters
 are called \textit{structured}.\index{structured population model}
Since living creatures reproduce from time to time, the evolution of these 
parameters has a ``jumping nature". 

Let $X\subset \mathbb R^d$ and let $x\in X$ be the vector of parameters 
that characterize any individual. 
We assume that the set $X$ has
non-empty interior and
has the boundary of measure zero.
We assume that  $x$ changes according to the equation
\begin{equation}
\label{s-p1}
x' (t) = b(x (t)),
\end{equation}
so the function $b$ and the set $X$ have such  properties that if $x_0\in X$, then there exists a unique solution of \eqref{s-p1} 
and $x(t)\in X$ for $t\ge 0$.

We assume that $\mu(x)$
is the death rate which includes both the real death and
the lost of cells during the process of division in cellular population models. 
Let $\psi(x)$  be the birth rate, that is 
 $\psi(x) \Delta t$ is the  probability that an individual 
with parameter $x$ has a  descendant
in time interval $[t,t+\Delta t]$. 
We assume that $\mu$ and $\psi$ are nonnegative and continuous functions
and that $\psi$ is bounded. 

We also assume that the transition probability function 
$\mathcal P(x,dy)$ shows the distribution
of parameters of a descendant
 at the birth. 
For example, if $x$ is the age then
$\mathcal P(x,\{0\})=1$;
if $x$ is the size of a cell, then 
 $\mathcal P(x,\{x/2\})=1$;
and in a phenotypic model the measure $\mathcal P(x,\cdot)$ 
is usually distributed near point $x$.

We are interested in finding the distribution
of the parameter $x$ at time $t$. This distribution
is described by a density function
$u(t,x)$; precisely 
\[
\int_A u(t,x)\,dx
\]
is the number of individuals (or biomass) with the parameter $x$ in the set  $A$.  
It should be underline that $u$ is not 
a density in the probabilistic sense because 
the integral of $u(t,x)$ over the whole phase space $X$
can be different from 1 and this integral can change with time.

First we  consider the case when $\mathcal P$ satisfies
\eqref{trana},
i.e. $\mathcal P$ corresponds to a stochastic operator $P$ 
on the space $L^1(X,\mathcal B(X),\ell)$.
Then $u$ satisfies an equation similar to \eqref{generator-A0}:
\begin{equation}
\label{sp-2}
\frac{\partial u}{\partial t}=-\nabla_x(bu)- \mu u+ P(\psi u)
\end{equation}
with initial condition $u(0,x)=u_0$. Moreover, if $P(t)u_0(x)=u(t,x)$,
then $\{P(t)\}_{t\ge 0}$ is a positive $C_0$-semigroup on 
$L^1(X,\mathcal B(X),\ell)$. 

The case when $\mathcal P$ does not satisfy \eqref{trana}
is generally more difficult. We restrict ourselves to an age-structured model,
i.e. a model with $\mathcal P(a,\{0\})=1$, $a$ denotes the age of an individual.
Then $u$ satisfies the following boundary-initial problem:
\begin{equation}
\label{McKm}
\begin{aligned}
&\frac{\partial u(t,a)}{\partial t}+\frac{\partial u(t,a)}{\partial a}= -\mu(a)u(t,a)\\
&u(t,0) =  \int_0^{\infty} u(t,a) \psi (a) \, d a,\quad u(0,a)=u_0.
\end{aligned}
\end{equation}
This problem also generates
 a positive $C_0$-semigroup.
 
The study of these models has a long history and there are many papers on them and their generalizations, see for example~\cite{BPR,DHT,GH,MD,RP,Webb}.
One of the key issues was  asynchronous exponential growth \eqref{c-d3}.
This property can be demonstrated 
using the spectral decomposition of quasi-compact positive semigroups (see Sec.~\ref{ss:expon-grow}).
This property can also be proved using theorems concerning stochastic semigroups. In this case the proof goes as follows.
Let $A$  be  an infinitesimal generator of
a positive semigroup 
$\{T(t)\}_{t\ge0}$ on $L^1(X,\mathcal B(X),\ell)$.  
We prove that there exist $\lambda\in\mathbb R$ and
 a bounded above and bounded away from zero 
function  $w$ such that 
$A^*w=\lambda w$. From this it follows that the semigroup
$\{P(t)\}_{t\ge0}$ given by $P(t)=e^{-\lambda t}T(t)$ is a
stochastic semigroup on the space $L^1(X,\mathcal B(X),m)$, where $m$ is a
measure  given by $m(B)=\int_B w(x)\,dx$. Next we prove 
asymptotic stability of the semigroup $\{P(t)\}_{t\ge0}$.
Since measures $\ell$ and $m$ are equivalent we have \eqref{c-d3}.

Dynamical systems with random jumps appear in many biological models.
Examples of such processes are \textit{ecological models with catastrophes}\index{ecological model with catastrophes} \cite{BGR,Gripenberg,HaT,Lande}.
The population size is described by a stochastic process $(\xi_t)_{t\ge 0}$, which between successive catastrophes
changes according to a differential equation
\begin{equation}
\label{kat1}
\xi_t'=b(\xi_t). 
\end{equation}
We assume that $\psi(x)$ is the intensity of occurrence of a catastrophe when the population size is $x$ 
and the transition probability function 
$\mathcal P(x,dy)$ shows the distribution of population size  after the catastrophe.
We assume that $\mathcal P$ corresponds to a stochastic operator $P$ 
on the space $L^1(\mathbb R_+,\mathcal B(\mathbb R_+),\ell)$.
If $u(t,x)$ is the density of the distribution of size of population then $u$ satisfies the equation
\begin{equation}
\label{e:kat2}
\frac{\partial u}{\partial t}=-\frac{\partial (bu)}{\partial x}-\psi u+P(\psi u)
\end{equation}
and generates a stochastic semigroup on 
$L^1(\mathbb R_+,\mathcal B(\mathbb R_+),\ell)$.

A similar model describes the \textit{dynamics of immune status}\index{dynamics of immune status} \cite{DGKT,PR-imun}.
The immune status is the concentration of specific antibodies, which appear after infection with a pathogen
and remain in serum, providing protection against future attacks by the same pathogen. 
Over time the number of  antibodies decreases until the next infection.  
During an infection, the immunity is boosted and then 
it is gradually waning, until it is boosted at the next infection.
Concentration of antibodies  between subsequent infections is a decreasing stochastic process $\xi_t$, which
satisfies Eq.~\eqref{kat1} with $b(x)<0$ for $x>0$.
The moments of infections are distributed according to a Poisson process with intensity $\lambda>0$.
If $x$ is the concentration of antibodies at the moment of infection, then $\varphi(x)>x$ is the concentration of antibodies 
just after clearance of infection. 
In this case the jump operator is the Frobenius--Perron operator $P_{\varphi}$ corresponding to $\varphi$ and $Af=-(bf)'+\lambda P_{\varphi}f-\lambda f$
is a generator of a stochastic semigroup on $L^1(\mathbb R_+,\mathcal B(\mathbb R_+),\ell)$.

\subsection{Semi-Markov kangaroo movement}
\label{ss:kanngaroo-movement}
Continuous-time Markov chains belong to a slightly larger class called
pure jump-type Markov processes.
A \textit{pure jump-type Markov process}\index{pure jump-type Markov process}
is a Markov process which remains constant between jumps.
Such processes are used to describe a movement in which an individual makes
 a series of jumps and the time of a single jump is neglected.
This process is also referred to as 
\textit{kangaroo process} or  \textit{kangaroo movement}~\cite{ODA}.\index{kangaroo movement}  

A pure jump-type homogeneous Markov process 
on a measurable space $(X,\Sigma)$ can be defined in the following way.
Let $\psi\colon X\to [0,\infty)$ be a given measurable function and let $\mathcal P$ be a
given transition probability function  on $X$.
Let $t_0=0$ and let $\xi_0$ be an $X$-valued random variable.
 For each $n\ge 1$ one can choose the $n$th
\textit{jump time} $t_n$ as a positive random variable
satisfying
\[
{\rm P}(t_n-t_{n-1}\le t|\xi_{t_{n-1}}=x)=1-e^{-\psi(x)t}, \quad t\ge 0,
\]
and we assume that $\xi_t=\xi_{t_{n-1}}$ for $t_{n-1}\le t<t_{n}$
and $\xi_{t_n}$ is a random $X$-valued variable  satisfying 
condition 
\[
{\rm P}(\xi_{t_n}\in B|\xi_{t_{n-1}}=x)=\mathcal P(x,B),\quad B\in\Sigma.
\]

If $m$ is a $\sigma$-finite measure on the space  $(X,\Sigma)$
and the transition probability function $\mathcal P$ satisfies condition~\eqref{trana},
then $\mathcal P$ corresponds to a jump operator $P$, which is a stochastic operator  on  
the space $L^1(X,\Sigma,m)$.
If  $\psi$ is a bounded function, then the operator
\begin{equation}
\label{op-infi-skok-kang-Mark}
Af(x)= P(\psi f)(x) - \psi(x)f(x)
\end{equation}
is a generator of uniformly continuous stochastic semigroup  
$\{P(t)\}_{t\ge 0}$ on $L^1(X,\Sigma,m)$. Indeed, if
 $\lambda\ge \sup\,\{\psi(x)\colon x\in X\}$
then we define the operator 
\[
\bar Pf(x)=\lambda^{-1}\big[P(\psi f)(x)+(\lambda -\psi(x))f(x)\big]
\]
and we represent the operator $A$ in the form 
$Af=\lambda \bar Pf-\lambda f$. The semigroup   $\{P(t)\}_{t\ge 0}$ takes the form
 (\ref{kang-wz-s2}): $P(t)f=e^{\lambda t(\bar P-I)}f$.

The time between jumps in the kangaroo movement described above  has an exponential distribution.
Let us now consider the movement of a kangaroo 
with the same jump operator $P$ as defined above, but
with other distributions of time between jumps. 
Such a process has no Markov property and it  belongs to the family of semi-Markov processes.

We assume that the \textit{holding time}\index{holding time}, i.e. the time  $a$ that the process spends at a given point $x$,  
has the probability density distribution $a\mapsto q(x,a)$.
We assume that $q$ is a measurable function and  we will need further technical  assumptions, which allow us to introduce some processes. 
In the case of previously introduced pure jump-type homogeneous Markov
process, $q$ is given by $q(x,a)=\psi(x)e^{-a\psi(x)}$.

Assume that the position of a jumping point is described by a stochastic process $(\xi_t)_{t\ge 0}$. The process $(\xi_t)_{t\ge 0}$ is 
generally non-Markovian, so it does not correspond to  a stochastic semigroup.
To overcome these obstacles we define  a Markov process $(\eta_t)_{t\ge 0}$ 
with the formula
\[
\eta_t=(\xi_{t_n},t-t_n) \quad\text{for $t_n\le t< t_{n+1}$,}
\]
where $t_0<t_1<t_2<t_3<\dots$ are jump times. 
If  $\xi_{t_n}=x$,
then the jump rate at time $t_n+a$ is given by
\begin{equation}
\label{def:p1}
p(x,a)=\lim_{\Delta t\downarrow 0}
\frac{\operatorname P(\tau\in [a,a+\Delta t] \mid \tau\ge a)}{\Delta t},
\end{equation}
where the random variable $\tau$ is the holding time. 
Then 
\begin{equation}
\label{def:p2}
p(x,a)=\frac{q(x,a)}{\int_a^{\infty} q(x,r)\,dr}
\end{equation}
and formally $p(x,a)=0$ when the last integral is zero.
The function $p$ plays the same role in this model as the function $\psi$ from the previous model except that it is a function of two variables $(x,a)$.

Here we assume that $t_0\le 0$ and start observing the process from the moment $t=0$, so the first jump occurs at a random time $t_1>0$.
Assume that the distribution of the random variable $\eta_0$ 
has density $u_0(x,a)$.
Then the random variable $\eta_t$
has some  density $u(t,x,a)$.
Since $p(x,a)$ is the rate of jump 
of the process $(\eta_t)_{t\ge 0}$ 
from $(x,a)$ to a point $(y,0)$, where $y$ is chosen according to the distribution $\mathcal P$, 
we have 
\begin{equation}
\label{boubd-cond}
u(t,x,0)=\int_0^{\infty}
\Big(P\big(p(\cdot,a)u(t,\cdot,a)\big)\Big)(x)\,da. 
\end{equation}
We will use the shortened notation $\mathfrak Pu(t,x)$
for the expression on the right-hand side of (\ref{boubd-cond}).
Thus, equation (\ref{boubd-cond}) takes the form 
$u(t,x,0)=\mathfrak Pu(t,x)$. 

Hence the function $u$ 
 satisfies the following initial-boundary problem:
\begin{align}
\label{eq1}
&\frac{\partial u}{\partial t}(t,x,a)
 +\frac{\partial u}{\partial a}(t,x,a)
 =-p(x,a)u(t,x,a),\\
&u(t,x,0)=\mathfrak Pu(t,x), 
\label{eq2}\\
&u(0,x,a)=u_0(x,a). 
\label{eq3}
\end{align}

The issue of when the resulting system of equations generates a stochastic semigroup is quite difficult and as we mentioned before we need additional assumptions. We assume that $X$ is an open  subset of $\mathbb R^d$
and the function $q$ satisfies the following conditions
\begin{enumerate}
\item[(A1)] $q\colon X\times [0,\infty)\to [0,\infty)$ is a  continuous function and 
for each $x\in X$ the function $a\mapsto q(x,a)$ is a probability density,
\item[(A2)] there exists $\varepsilon>0$ such that
$\int_0^{\varepsilon} q(x,a)\,da<1-\varepsilon $ for all $x\in X$.
\end{enumerate}

\def\ga{\overline a}
We define $\ga(x)\le \infty$ to be the minimum number such that
\[
\int_0^{\ga(x)}q(x,a)\,da=1\quad\textrm{for $x\in X$}.
\]
Let
\begin{equation}
\label{eq:Y}
Y=\{(x,a)\colon    x\in X,\,\,\, 0\le a<\ga(x)\}.
 \end{equation}
The process $(\eta_t)_{t\ge 0}$ has values in  $Y$. 

Let  $E=L^1(Y,\mathcal B(Y),m)$, where $\mathcal B(Y)$ is the $\sigma$-algebra of Borel subsets of $Y$ and $m$ is the Lebesgue measure.
We replace problem (\ref{eq1})--(\ref{eq3}) by an evolution equation $u'(t)=Au(t)$ on the space $E$ with operator $A$ of the form 
\begin{equation}
\label{eq:A}
 Af=-\frac{\partial f}{\partial a}-pf.
\end{equation}
and  with domain
\begin{equation}
\label{domain-A}
\mathfrak D(A)=\Big\{f\in E\colon\,\, \frac{\partial f}{\partial a}\in E,
\,\,pf\in E,\,\, f(x,0)=\mathfrak Pf(x) \Big\}.
 \end{equation}

The operator $A$ generates a stochastic semigroup $\{P(t)\}_{t\ge0}$
on the space~$E$. We omit the proof of this result here.
The interested reader is referred to publication \cite{PR-cell-2022} where the proof of a similar result is given. This proof is based on  
 perturbation method related to
operators with boundary conditions developed in \cite{greiner}
and an extension of this method to unbounded perturbations in $L^1$ space
in \cite{GMTK}.

\subsection{Cell cycle model}
\label{ss:cykl-kom}
The cell cycle is a series of events that take place in a cell leading to its replication. It is regulated by a
complex network of protein interactions~\cite{Morgan}.
One of the most important questions concerning the cell cycle is what factor determines when a cell divides.
Such a factor may be the size of the cell and therefore different size-structured models are considered.
Such models can be experimentally verified and allow us  to understand
processes inside single cells.

We assume that the population  grows in steady-state conditions and neglected mortality of cells. 
Cells can be described by their age $a$ and size $x$ alone and reproduction occurs by fission into two equal parts. 
We assume that a cell grows with velocity $g(x)$, i.e. $x'(t)=g(x(t))$.

First, we assume that also the rate of division $\lambda$ depends only on a cell size, then the evolution of size distribution of cells is described by
the model considered in Sec.~\ref{ss:structured}, where 
$\mathcal P(x,\{x/2\})=1$, $\mu\equiv 0$ and $\psi(x)=2\lambda(x)$. 
In the definition of $\psi$ we have a factor $2$ 
because division results in two cells.
Moreover, $Pf(x)=2f(2x)$ is the stochastic operator $P$ corresponding to $\mathcal P$. Thus the density function $u(t,x)$ satisfies a version of
Eq.~\eqref{sp-2}:
\begin{equation}
\label{sp-cc}
\frac{\partial u(t,x)}{\partial t}=
-\frac{\partial \big(g(x)u(t,x)\big)}{\partial x}
+ 4\lambda(2x)u(t,2x).
\end{equation}

Although models of type \eqref{sp-cc} are commonly used, they have 
one disadvantage in that it is difficult to determine
the function $\lambda$  experimentally.
In~\cite{PR-cell-2022} we propose a model that is based on two factors: the initial size of a cell $x_b$ and  
the length of the cell cycle. Both factors are not difficult to measure 
and there are quite a few results and models in the literature based on experiments on the distribution of cell cycle length and its dependence on initial cell size.
However, a model based on the length of the cell cycle does not lead to an equation for $u(t,x)$, but using its properties we can investigate the long-time behaviour of $u(t,x)$.

We now present this model and give its main properties.
We assume that sizes of cells and cell cycle durations are bounded
above and bounded away from zero. 
Let  $\dxb$ and $\gxb$   be the minimum and maximum sizes of newborn cells.
We also assume that cells age with unitary velocity and grow  with a velocity $g(x)$,
i.e. if a cell has the initial size $x_b$, then the size at age $a$ satisfies the equation
\begin{equation}
\label{grow}
x'(a)=g(x(a)),\quad x(0)=x_b.
 \end{equation}
We denote by $\pi_ax_b$ the  solution of (\ref{grow}).
If the initial size of the mother cell is $x_b$ and $\tau=a$ is the length of its cell cycle, then the initial size of the daughter cell is $S_{a}(x_b)=\tfrac12\pi_ax_b$.

The length $\tau$ of the cell  cycle is a random variable which depends on the initial cell size $x_b$; has values in some interval 
$[\da(x_b),\ga(x_b)]$;  and has the probability density distribution $q(x_b,\cdot)$, i.e.
the integral $\int_0^{A} q(x_b,a)\, da$ is the probability that $\tau\le A$.  
According to the definition of $q$, if a cell has the initial size $x_b$, then $\Phi(x_b,a)=\int_a^{\infty} q(x_b,r)\,dr$ is its \textit{survival function},
i.e. $\Phi(x_b,a)$ is the probability that a cell will not split before age $a$.

The functions $g$ and $q$ satisfy the following assumptions:
\begin{enumerate}
\item [(A1)] $g\colon [\dxb,2\gxb]\to (0,\infty)$
 is a $C^1$-function,

\item [(A2)] $q\colon [\dxb,\gxb]\times [0,\infty)\to [0,\infty)$ is a continuous function and for each $x_b$ the function $a\mapsto q(x_b,a)$ is a probability density,

\item [(A3)]  $0<\da(x_b)<\ga(x_b)<\infty$,  $q(x_b,a)>0$ if $a\in  (\da(x_b),\ga(x_b))$, 
and  $q(x_b,a)=0$ if $a\notin (\da(x_b),\ga(x_b))$ for each $x_b\in [\dxb,\gxb]$, 

\item [(A4)] $x_b\mapsto \da(x_b)$ and $x_b\mapsto \ga(x_b)$ are continuous functions,

\item [(A5)]  $S_{\da(x_b)}(x_b)\ge \dxb$ and $S_{\ga(x_b)}(x_b)\le \gxb$ for each $x_b\in [\dxb,\gxb]$,

\item [(A6)] $S_{\da(x_b)}(x_b)<x_b<S_{\ga(x_b)}(x_b)$ for each $x_b\in (\dxb,\gxb)$. 
\end{enumerate}
Assumption (A6) means that the initial daughter cell size is  distributed in some interval which contains the mother initial size.

Let the function $p$ be given by  \eqref{def:p1}.
Then $p$ satisfies \eqref{def:p2} and
$q=p\Phi$. Let $P_a$ be the Frobenius--Perron operator corresponding to the transformation $S_a$.
Denote by $u(t,x_b,a)$ the number of individuals in a population  having initial size $x_b$ and age $a$ at time $t$.
Then,  according to our assumptions concerning the model, 
$p(x_b,a)u(t,x_b,a)$ is the number of cells of initial size $x_b$ and age $a$ which split in a unit time interval. It means that 
\[
\mathfrak Pu(t,x_b)= 2\int_0^{\infty}\Big(P_a\big(p(\cdot,a)u(t,\cdot,a)\big)\Big)(x_b)\,da 
\]
 is the number of new born cells in this time interval. 
Thus the function $u$ satisfies 
the initial-boundary problem \eqref{eq1}--\eqref{eq3}.

As in Sec.~\ref{ss:kanngaroo-movement}
we replace problem (\ref{eq1})--(\ref{eq3}) by an evolution equation $u'(t)=Au(t)$ on the space 
$E=L^1(Y,\mathcal B(Y),m)$
 with generator $A$ given by \eqref{eq:A} and  \eqref{domain-A}.
 Here the set $Y$ is given by \eqref{eq:Y}.
The operator $A$ generates a positive $C_0$-semigroup $\{T(t)\}_{t\ge0}$
on the space~$E$.

In order to formulate the main result of \cite{PR-cell-2022}
we need an additional assumption:
\begin{enumerate}
\item[(A7)] there exists $x\in (\dxb,\gxb)$ such that $g(2x)\ne 2g(x)$.
\end{enumerate}

\begin{theorem}
\label{AEG-cc}
If conditions {\rm (A1)}--{\rm (A7)} are fulfilled, then the semigroup $\{T(t)\}_{t\ge0}$
has asynchronous  exponential growth, i.e.
there exist $\lambda\in\mathbb R$,
a  density $f^*\in E$ and a positive function 
$v\in L^{\infty}(Y,\mathcal B(Y),m)$
such that 
\begin{equation}
\label{AEG-cc-w}
\lim_{t\to\infty}e^{-\lambda t}T(t)f=
f^* \iint\limits_Y f(x_b,a)v(x_b,a)\,dx_b\,da  \quad\textrm{for $f\in E$}.
\end{equation}
\end{theorem}
The proof of Theorem~\ref{AEG-cc} is based on 
the method presented in Sec.~\ref{ss:structured} of reducing the problem to the asymptotic stability of some stochastic semigroup,   
and then applying Theorem~\ref{asym-th2}.

Now we consider a size-age-version of our model. 
Let  $w(t,x,a)$ be the number of cells   having  size $x$ and age $a$ at time $t$. Then 
 \[
u(t,x,a)=\frac{\partial (\pi_ax)}{\partial x} w(t,\pi_a x,a)
 =\frac{g(\pi_ax)}{g(x)} w(t,\pi_a x,a).
\]
From the formula \eqref{AEG-cc-w} it follows that
\[
e^{-\lambda t} w(t,x,a)\to h^*(x,a)
   \iint\limits_Z v(\pi_{-a}y,a)w(0,y,a)\,dy\,da
\]
in $E$ as $t\to\infty$, where $h^*(x,a)=\dfrac{g(\pi_{-a}x)}{g(x)}f^*(\pi_{-a}x,a)$
and 
\[
Z=\{(\pi_ax,a)\colon    x\in X,\,\,\, 0\le a<\ga(x)\}.
\]
In  particular, if $\bar w(t,x)= \int_0^{\infty}w(t,x,a)\,da$,
then 
\[
e^{-\lambda t} \bar w(t,x)\to \bar h^*(x)
   \iint\limits_Z v(\pi_{-a}y,a)w(0,y,a)\,dy\,da,\quad 
\bar h^*(x)=\int_0^{\infty}h^*(x,a)\,da.
\]
Thus, we have obtained an asymptotic distribution of cell sizes, even though the function $\bar w$ is not determined by a semigroup of operators.

\begin{acknowledgement}
This research was partially supported by
the  National Science Centre (Poland)
Grant No. 2017/27/B/ST1/00100.
\end{acknowledgement}

\printindex
\end{document}